\documentclass [11pt,oneside]{article}

\reversemarginpar \textwidth=14cm \textheight=19cm

\usepackage{amssymb} \usepackage{amsfonts} \usepackage{amsmath}
\usepackage{amsthm} \usepackage{epsfig}

\newcommand{\mettifig}[1]{\epsfig{file=#1}}
\renewcommand{\leq}{\leqslant}
\renewcommand{\geq}{\geqslant}
\usepackage{graphicx}

\newtheorem{lemma}{Lemma}[section] 
\newtheorem{teo}[lemma]{Theorem}
\newtheorem{rem}[lemma]{Remark} 
\newtheorem{prop}[lemma]{Proposition}
\newtheorem{cor}[lemma]{Corollary}

\newtheorem{quest}[lemma]{Question}

\newcommand{\matN}{\ensuremath {\mathbb{N}}}
\newcommand{\matR} {\ensuremath {\mathbb{R}}}
\newcommand{\matQ} {\ensuremath {\mathbb{Q}}}
\newcommand{\matZ} {\ensuremath {\mathbb{Z}}}
\newcommand{\matC} {\ensuremath {\mathbb{C}}}

\newcommand{\matH} {\ensuremath {\mathbb{H}}}

\newcommand{\calI} {\ensuremath {\mathcal{I}}}

\newcommand{\calM} {\ensuremath {\mathcal{M}}}

\newcommand{\calP} {\ensuremath {\mathcal{P}}}
\newcommand{\calT} {\ensuremath {\mathcal{T}}}

\newcommand{\calO}{\ensuremath {\mathcal{O}}}
\newcommand{\calG}{\ensuremath {\mathcal{G}}}

\newcommand{\TV}{{\rm TV}}

\newcommand{\nota} [1] {\caption{\footnotesize{#1}}}

\newfont{\Got}{eufm10 scaled 1200}

\font\titsc=cmcsc10 scaled 1200

\newcommand{\finedimo}{{\hfill\hbox{$\square$}\vspace{2pt}}}

\newcommand{\diag}{\SelectTips{cm}{} \xymatrix@1}

\author{Roberto \titsc{Frigerio} \and Bruno \titsc{Martelli}
\and Carlo \titsc{Petronio}}

\title{Dehn filling of cusped hyperbolic 3-manifolds \\
with geodesic boundary}

\newcommand{\tetra}{\Delta}
\newcommand{\slopedist}{\Delta}
\newcommand{\hypbound}{\Delta^{\rm hyp}}
\newcommand{\largebound}{\Delta^{\rm large}}
\newcommand{\negbound}{\Delta^{\rm neg}}

\begin{document}

\maketitle

\begin{abstract}
\noindent
	We define for each $g\geqslant 2$ and $k\geqslant 0$
	a set $\calM_{g,k}$ of orientable hyperbolic 3-manifolds with $k$ toric cusps and a
	connected totally geodesic boundary of genus $g$. Manifolds in $\calM_{g,k}$
	have Matveev complexity $g+k$ and Heegaard genus $g+1$, and their homology, volume, and Turaev-Viro 
	invariants depend
	only on $g$ and $k$. In addition, they do not contain closed essential surfaces.
	The cardinality of $\calM_{g,k}$ for a fixed $k$ has growth type $g^g$.
	
	We completely describe  the non-hyperbolic Dehn fillings of 
each $M$ in $\calM_{g,k}$, showing that, on any cusp of any hyperbolic 
manifold obtained by partially filling $M$, there are precisely $6$ non-hyperbolic Dehn
	fillings: three contain essential discs,
	and the other three contain essential annuli. This gives an infinite class
	of large hyperbolic manifolds (in the sense of Wu) with $\partial$-reducible
	and annular Dehn fillings having distance 2, and allows us to prove that the
	corresponding upper bound found by Wu is sharp.
	If $M$ has one cusp only, the three $\partial$-reducible fillings are handlebodies.
	
  \vspace{4pt}

\end{abstract}

\noindent

\section{Definition and statements} \label{statements:section}
In this paper we introduce certain classes $\calM_{g,k}$ of 
compact 3-manifolds, we determine many topological and geometric
invariants of the elements of $\calM_{g,k}$, and we analyze their Dehn
fillings, answering in particular a question raised by Wu~\cite{Wu}
on the distance between non-hyperbolic fillings of a large 3-manifold.
We also show that $\#\calM_{g,k}$ grows very fast as $g$ goes to infinity.

\paragraph{Definition of $\calM_{g,k}$} 
All the manifolds considered in this paper
will be viewed up to homeomorphism, and will be connected and orientable by default.
Let $\tetra$ denote the standard tetrahedron, and let $\dot{\tetra}$ be
$\tetra$ with its vertices removed.
An \emph{ideal triangulation} of a compact 3-manifold $M$ with boundary
is a realization of the interior of $M$ as a gluing of a finite
number of copies of $\dot{\tetra}$, induced by a simplicial face-pairing  
of the corresponding $\tetra$'s.
Let $\Sigma_g$ be the closed orientable surface of genus $g$. The following
result is proved in Section~\ref{triangulations:section} and 
motivates our definition of $\calM_{g,k}$.
\begin{prop} \label{first:prop}
An ideal triangulation of a manifold whose boundary is the union of $\Sigma_g$
and $k$ tori
contains  at least $g+k$ tetrahedra.
\end{prop}
\noindent
We then define $\calM_{g,k}$ for all $g\geqslant 2, k\geqslant 0$ as follows:
\begin{eqnarray*}
\calM_{g,k} & = & 
\big\{ {\rm compact\ orientable\ manifolds}\ M\ {\rm having\ an\ ideal\ triangulation}\\ 
& & \quad {\rm with\ } g+k\ {\rm tetrahedra,\ and\ }
\partial M = \Sigma_g\sqcup\big(\mathop{\sqcup}\limits_{i=1}^k T_i\big)\ 
{\rm with}\ T_i\cong\Sigma_1\big\}.
\end{eqnarray*}
The sets $\calM_g = \calM_{g,0}$ were studied in~\cite{FriMaPe1}.

\paragraph{Geometric and topological invariants}
We now describe the main properties of the manifolds in $\calM_{g,k}$, starting from
a quick general review of the invariants that we can compute.

We first recall that a 
surface in a compact 3-manifold $M$ is
\emph{essential} if it is
properly embedded, connected, and either
a reducing sphere, or a boundary-reducing disc, or an
incompressible and $\partial$-incompressible surface not parallel to
the boundary.

We say that a compact $3$-manifold $M$ is
\emph{hyperbolic} if, after removing the boundary tori, we get
a complete hyperbolic $3$-manifold with finite volume and geodesic
boundary. If $\partial M\neq\emptyset$, Thurston's
geometrization theorem implies that hyperbolicity is
equivalent to the condition that $M$ does
not contain any essential surface with non-negative Euler characteristic.
We recall that Kojima has proved in~\cite{Kojima} that every hyperbolic manifold with non-empty
geodesic boundary admits a \emph{canonical decomposition} into geometric polyhedra.

For any compact $3$-manifold $M$, an $\matN$-valued invariant $c(M)$
was defined by Matveev
in~\cite{Mat} and called the \emph{complexity} of $M$. Matveev also proved that,
when $M$ is hyperbolic, $c(M)$ equals the minimal number of
tetrahedra in an ideal triangulation of $M$. 

If $M$ is a compact 3-manifold with $\partial M=\partial_0M\sqcup\partial_1M$,
one can define the \emph{Heegaard genus} of $(M,\partial_0M,\partial_1M)$ as the minimal
genus of a surface that splits $M$ as $C_0\sqcup C_1$, where $C_i$
is obtained by attaching $1$-handles on one side of a collar of $\partial_iM$.

For any compact 3-manifold $M$ and integer $r\geqslant 2$,
after fixing in $\matC$ a primitive $2r$-th root of unity,
a real-valued invariant $\TV_r(M)$
for compact 3-manifolds with boundary was defined by Turaev and Viro in~\cite{TuVi}.

The following theorem will be proved in Sections~\ref{triangulations:section}
and~\ref{spines:section}.

\begin{teo} \label{main:teo}
Let $M\in\calM_{g,k}$. The following holds:
\begin{enumerate}
\item \label{main:hyperbolic:item} 
$M$ is hyperbolic, and its volume depends only on $g$ and $k$;
\item \label{main:canonical:item} 
$M$ has a unique ideal triangulation with $g+k$ tetrahedra, which gives
the canonical Kojima decomposition of $M$;
\item \label{main:incompressible:item}
Every essential surface in $M$ has 
non-empty boundary which intersects $\Sigma_g$;
\item \label{main:complexity:item}
$M$ has complexity $g+k$;
\item \label{main:genus:item}
The Heegaard genus of $\big(M,\Sigma_g,\mathop{\sqcup}\limits_{i=1}^kT_i\big)$ is $g+1$;
\item \label{main:homology:item}
$H_1(M;\matZ)=\matZ^{g+k}$;
\item \label{main:turaev:viro:item} 
The Turaev-Viro invariant $\TV_r(M)$ of $M$ depends only on $r$, $g$ and $k$.
\end{enumerate}
\end{teo}

\paragraph{Growth of $\calM_{g,k}$}
We begin with the following fact, established in Section~\ref{estimates:section}:

\begin{prop}\label{non-empty:prop} 
$\calM_{g,k}$ is non-empty precisely for $g>k$ or $g=k$ and $g$ even.
\end{prop}

We have computed the number of elements of $\calM_{g,k}$ for some
small $g$ and $k$ with the aid of a computer. 
Our results are summarized in the next table, where we also take
from~\cite{FriMaPe2} 
the number
of all hyperbolic manifolds with non-empty geodesic boundary having 
a certain complexity $c$:
\begin{center}
\makebox{
\begin{tabular}{c||c|c|c||c}
\phantom{\Big|}       & $\calM_{c,0}$ & $\calM_{c-1,1}$  & $\calM_{c-2,2}$  & All hyperbolic manifolds  \\ \hline\hline
\phantom{\Big|} $c=2$ & $8$           & $\emptyset$      & $\emptyset$      & $8$ \\ \hline
\phantom{\Big|} $c=3$ & $74$          & $1$              & $\emptyset$      & $151$ \\ \hline
\phantom{\Big|} $c=4$ & $2340$        & $12$             & $1$              & $5033$ \\ \hline
\phantom{\Big|} $c=5$ & ?             & $416$              & $1$              & ? \\ \hline
\phantom{\Big|} $c=6$ & ?             & $17900$                & $51$              & ? \\ 
\end{tabular}
}
\end{center}

We now say that a numerical sequence $\big(a_n\big)_{n=1}^\infty$ has 
\emph{growth type $n^n$} if there
exist constants $C>c>0$ such that $n^{c\cdot n}
< a_n < n^{C\cdot n}$ for $n\gg 0$. In Section~\ref{estimates:section}
we will prove the following:

\begin{teo} \label{growth:teo} 
For any fixed $k$ the sequence $\big(\#\calM_{g,k}\big)_{g=2}^\infty$ has growth type $g^g$.
\end{teo}

This result and an easy upper bound also established in Section~\ref{estimates:section}
readily imply the following:
\begin{cor} The number of hyperbolic $3$-manifolds of complexity
$c$ has growth type $c^c$.
\end{cor}

\paragraph{Dehn fillings}
Recall that a \emph{slope} in a torus $T$ is
an isotopy class of simple closed
essential curves, and that, after choosing a $\matZ$-basis of $H_1(T;\matZ)$, a slope
is represented by a number in $\matQ\cup\{\infty\}$.
If $M$ is a manifold with $k$ boundary tori, 
and $\alpha_1,\ldots,\alpha_h$ are slopes in some $h\leqslant k$
of these tori, we denote by
$M(\alpha_1,\ldots,\alpha_h)$ the result of Dehn-filling these $h$ tori 
along $\alpha_1,\ldots,\alpha_h$.

For $g\geq 2$ we now denote by $H_g$ the handlebody
of genus $g$ and for $0\leq k\leq g$ we introduce another manifold
$H_{g,k}$. We do this noting that $H_g$ can be viewed as the 
$\partial$-connected sum of $g$ solid tori, and defining
$H_{g,k}$ to be $H_g$ minus open tubes around the cores of $k$ of these
solid tori. So $H_{g,k}$ is obtained from $H_g$ by drilling out $k$ tunnels 
along $k$ different 1-handles.  Of course $H_{g,k}$ 
is well-defined and $H_{g,0}=H_g$. Moreover, $H_{g,k}$ is not hyperbolic because it is
$\partial$-reducible.  

The next result is proved in Section~\ref{spines:section}.

\begin{teo} \label{main:Dehn:teo}
Let $M\in\calM_{g,k}$, with 
$\partial M = \Sigma_g\sqcup\big(\mathop{\sqcup}\limits_{i=1}^k T_i\big)$.
There exists a $\matZ$-basis
of $H_1\big(\mathop{\sqcup}\limits_{i=1}^k T_i;\matZ\big)$ such that
$N=M(\alpha_1,\ldots,\alpha_h)$ is as follows:
\begin{itemize}
\item If $\alpha_i\in\{0,1,\infty\}$ for some $i$ then
$N = H_{g,k-h}$, so it is not hyperbolic;
\item If $\alpha_i\in\{-1,1/2,2\}$ for some $i$ then
$N$ contains a M\"obius strip or non-separating
annulus $R$ with $\partial R\subset\Sigma_g$, and cutting $N$ along $R$ one gets $H_{g,k-h}$;
also in this case $N$ is not hyperbolic; 
\item If $\alpha_i\not\in\{-1,0,1/2,1,2,\infty\}$ for all $i$ then
$N$ is hyperbolic and, denoting by $T_{j_1},\ldots,T_{j_{k-h}}$
the non-filled tori, the Heegaard genus of 
$$\big(N,\Sigma_g,T_{j_1}\sqcup\ldots\sqcup T_{j_{k-h}}\big)$$
is $g+1$;
\item If $\alpha_i\in\{-2,-1/2,1/3,2/3,3/2,3\}$ for all $i$ then
$N$ belongs to $\calM_{g,k-h}$.
\end{itemize}
Moreover every essential surface in $N$ has non-empty boundary intersecting $\Sigma_g$.
\end{teo}

If $\alpha,\alpha'$ are two slopes on a torus, we denote now by
$\slopedist(\alpha,\alpha')$ their \emph{distance}, that is their
geometric intersection number. We recall that, once
a homology basis is fixed, the set $\matQ\cup\{\infty\}$ of slopes
can be viewed as a subset of $\partial\matH^2$, where $\matH^2$ is the
half-space model of hyperbolic plane. Connecting
the pairs of slopes $\alpha,\alpha'$ such that $\slopedist(\alpha,\alpha')=1$
one gets a tessellation of $\matH^2$ by ideal triangles.
A combinatorial (but geometrically incorrect)
picture of this tessellation is shown in
Fig.~\ref{tessellation:fig} in the disc model of $\matH^2$.  
If $\slopedist(\alpha,\alpha')>1$, then $\slopedist(\alpha,\alpha')$
is computed from Fig.~\ref{tessellation:fig} as $1$ plus the number of
lines met by the line which joins $\alpha$ and $\alpha'$.
\begin{figure}[t,scale=.7]
\begin{center}
\input{tesselation.pstex_t}
\nota{The Farey tessellation of $\matH^2$.}\label{tessellation:fig}
\end{center}
\end{figure}
The theorem just stated gives the following:

\begin{cor} \label{main:cor}
For any $g\geqslant 2$ there exist infinitely many hyperbolic 
manifolds $N$ with $\partial N=\Sigma_g\sqcup\Sigma_1$ and with $6$ slopes
$\alpha^1, \ldots, \alpha^6$ on $\Sigma_1$, such that $N(\alpha^i) = H_g$
for $i\in\{1,2,3\}$ and $N(\alpha^i)$ is
annular for $i\in\{4,5,6\}$. We have $\slopedist(\alpha^i, \alpha^{i+3})=2$
for $i=\{1,2,3\}$, and $\slopedist(\alpha^j,\alpha^{j'}) = 3$ for
$j,j'\in\{4,5,6\}, j\neq j'$.
\end{cor}

\begin{proof}
Take $N = M(\alpha)$ where $M \in \calM_{g,2}$ and $\alpha$ varies
in $\matQ\setminus\{-1,0,1/2,1,2\}$.
This gives infinitely many manifolds because the volume grows
to $\mathrm{vol}(M)$ as 
$\slopedist(\alpha,0)$ tends to infinity. Now let 
$(\alpha^1,\ldots,\alpha^6) = (0,1,\infty,2,-1,1/2)$, so
$N(\alpha^i) = H_g$ for $i\in\{1,2,3\}$. For  
$i\in\{4,5,6\}$ the manifold 
$N(\alpha^i)$ is not hyperbolic, 
so by Thurston's hyperbolization
theorem and the last assertion of Theorem~\ref{main:Dehn:teo}
it contains
either an essential disc or an essential annulus. Since a
hyperbolic manifold admits at most $3$ boundary-reducible
fillings~\cite{Wu2}, we conclude that $N(\alpha_i)$ is annular
for $i\in\{4,5,6\}$.
\end{proof}

This corollary leads to infinitely many examples of knots in a handlebody $H_g$ having 
non-meridinal fillings that give back $H_g$ (more precisely, two such
fillings) for all $g\geqslant 2$.
For $g=1$, \emph{i.e.} for the solid torus $H_1$, 
knots in $H_1$ with non-meridinal fillings giving $H_1$ 
were shown to be 1-bridge~\cite{Ga2}
and then classified by Gabai~\cite{Ga} and Berge~\cite{Be}.
(See Section~\ref{spines:section} for a definition of $1$-bridge knot). 
In particular, there is a unique 
knot in $H_1$ with two non-meridinal fillings giving $H_1$. 
A knot in $H_2$ with one non-meridinal filling giving $H_2$ was
also shown in~\cite{Ga}, and in the same paper 
other examples in $H_g$ for any $g$
were attributed to Berge, together with the following question: if
$K$ is a knot in a $\partial$-reducible manifold $M$ with
$\partial$-irreducible exterior and a $\partial$-reducible non-meridinal surgery, is
$K$ boundary-parallel or $1$-bridge in $M$? 
The answer is ``yes'' for all the hyperbolic examples one can
construct from $\calM_{g,k}$, as we will show
in Section~\ref{spines:section}:

\begin{prop} \label{1-bridge:prop}
Every knot in $H_g$ whose complement is a hyperbolic manifold obtained from $k-1$ Dehn
fillings on a manifold in $\calM_{g,k}$ is 
a $1$-bridge knot.
\end{prop}

Let $M$ be a hyperbolic $3$-manifold, and 
$T_i\subset\partial M$ be a chosen boundary torus. A slope
$\alpha$ on $T_i$ is called \emph{exceptional} if $M(\alpha)$ is not hyperbolic.
Assuming $M$ has either some other cusp or non-empty geodesic
boundary, $\alpha$ is exceptional if and only if $M(\alpha)$ contains an
essential sphere, disc, annulus, or torus. In these cases we say that $\alpha$ is 
respectively of type $S$, $D$, $A$, or $T$.
For $X_1,X_2\in\{S,D,A,T\}$, we define a number 
$\hypbound(X_1, X_2)$ as the maximal distance of two
slopes of type $X_1$ and $X_2$ in a boundary torus of a hyperbolic manifold~\cite{Go}.
Wu pointed out~\cite{Wu} that in most cases $\hypbound$
is considerably lower when considering only \emph{large} manifolds,
\emph{i.e.}~manifolds with
$H_2(M,\partial M\setminus T_i)\neq\{0\}$.
He thus defined $\largebound(X_1,X_2)$ as the maximal distance
of two slopes of type $X_1$ and $X_2$ in a boundary torus of a
\emph{large} hyperbolic manifold. Among other inequalities, he proved that
$1\leqslant\largebound(D,A)\leqslant 2$.
Since every hyperbolic manifold $M$ with $\chi(M)<0$ is large,
Corollary~\ref{main:cor} implies that $\largebound(D,A) = 2$. 
This result leaves $\largebound(T,T)$ as the only unknown value
for $\hypbound$ and $\largebound$, as shown in the
following tables (which are taken from~\cite{Go} with the insertion
of $\largebound(D,A) = 2$).

\begin{center}
\makebox{
\begin{tabular}{@{} c||c|c|c|c}
\phantom{\Big|} $\hypbound$ & $S$ & $D$ & $A$ & $T$ \\ \hline\hline
\phantom{\Big|} $S$      & $1$ & $0$ & $2$ & $3$ \\ \hline
\phantom{\Big|} $D$      &     & $1$ & $2$ & $2$ \\ \hline
\phantom{\Big|} $A$      &     &     & $5$ & $5$ \\ \hline
\phantom{\Big|} $T$      &     &     &     & $8$ \\
\end{tabular}
}
\hspace{1 cm}
\makebox{
\begin{tabular}{@{} c||c|c|c|c}
\phantom{\Big|} $\largebound$ & $S$ & $D$ & $A$ & $T$ \\ \hline\hline
\phantom{\Big|} $S$        & $0$ & $0$ & $1$ & $1$ \\ \hline
\phantom{\Big|} $D$        &     & $1$ & $2$ & $1$ \\ \hline
\phantom{\Big|} $A$        &     &     & $4$ & $4$ \\ \hline
\phantom{\Big|} $T$        &     &     &     & $4$--$5$ \\
\end{tabular}
}
\end{center}

We now define $\negbound(X_1,X_2)$ as the maximal 
distance of two slopes of type $X_1$ and $X_2$ in a boundary torus
of some hyperbolic manifold $M$ with $\chi(M)<0$. Of course we have 
$$\negbound(X_1,X_2)\leqslant
\largebound(X_1,X_2).$$
Gordon and Wu proved in~\cite{GoWu-AA} that if two
slopes of type $A$ in a boundary torus of $M$ have distance greater than $3$
then $M$ is the complement of some link in $S^3$. Hence 
$\negbound(D,A)\leqslant 3$. This estimate and 
Corollary~\ref{main:cor} give the following results for $\negbound$:
\begin{center}
\begin{tabular}{@{} c||c|c}
\phantom{\Big|} $\negbound$    & $D$ & $A$ \\ \hline\hline
\phantom{\Big|} $D$            & $1$ & $2$ \\ \hline
\phantom{\Big|} $A$            &     & $3$ \\ 
\end{tabular}
\end{center}
More precisely,
the values of $\negbound(X_1,X_2)$ for
$X_1, X_2\in\{D,A\}$ are  realized by every manifold $M\in\calM_{g,k}$ with
$k\geqslant 1$, and by every hyperbolic manifold obtained from 
such an $M$ by filling some (but not all) boundary components.

\begin{quest}
What is $\negbound(X_1,X_2)$ when $X_1$ is
$S$ or $T$ ?
\end{quest}

\begin{figure}
    \begin{center}
    \mettifig{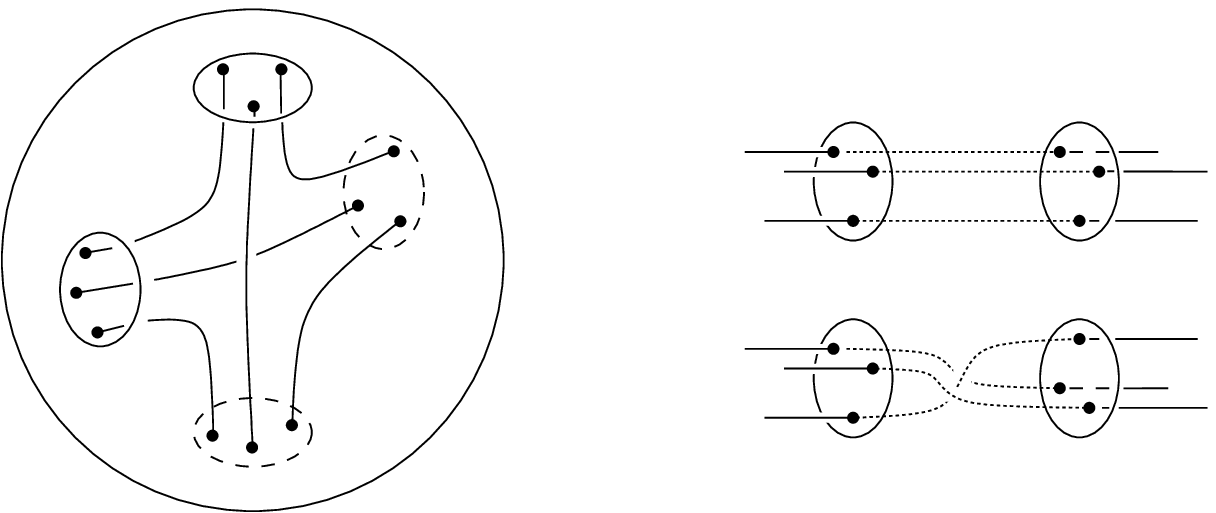, width = 9 cm}
    \nota{A tangle in a ball (left); a positive and a negative gluing of tangles (right).} 
    \label{tangle:fig}
    \end{center}
    \end{figure}

\paragraph{Knots giving $\calM_{g,1}$} It turns
out that the elements of $\calM_{g,1}$ are
knot exteriors in $H_g$, and the knots can be exhibited explicitly,
as we now explain. Consider a ball $B$ as in Fig.~\ref{tangle:fig}-left
with the tangle $\tau\subset B$ as shown, the 12 ends of $\tau$ being arranged 
in four groups of three, each group contained in a disc. Now take $g-1$ copies of 
$(B,\tau)$ and glue together the $4(g-1)$ discs, matching the ends
of the $\tau$'s. Each gluing should be of one of the two types
suggested in Fig.~\ref{tangle:fig}-right. A gluing 
as in the top part of Fig.~\ref{tangle:fig}-right
will be called \emph{positive}, one as in the bottom part will be called
\emph{negative}. The result of the $2(g-1)$ gluings is a link in $H_g$,
and one readily sees that if all the gluings are positive then
the link is parallel to $\partial H_g$. The 
following will be proved at the end of Section~\ref{spines:section}:

\begin{prop} \label{tangle:prop}
If a knot $K$ in $H_g$ is realized from $g-1$ copies of $(B,\tau)$ with
$2g-3$ positive gluings and one negative gluing then the exterior of $K$
belongs to $\calM_{g,1}$. Every manifold in
$\calM_{g,1}$ arises like this.
\end{prop}

\section{Triangulations and hyperbolicity} \label{triangulations:section}
In this section we discuss some basic properties of the manifolds in 
$\calM_{g,k}$. We describe in Proposition~\ref{forma:tria:prop} the properties of a 
triangulation of a compact manifold $M$ with
$\partial M=\Sigma_g\sqcup\big(\mathop{\sqcup}\limits_{i=1}^k T_i\big)$, 
showing in particular that such a triangulation has at least
$g+k$ tetrahedra. This result proves Proposition~\ref{first:prop} and,
together with hyperbolicity of the manifolds in
$\calM_{g,k}$, easily implies 
point~(\ref{main:complexity:item}) of Theorem~\ref{main:teo}.
We then prove all other points of Theorem~\ref{main:teo},
except point~(\ref{main:homology:item}), which is deferred to Section~\ref{spines:section}.
Namely, solving the
hyperbolicity equations~\cite{FriPe} we prove point~(\ref{main:hyperbolic:item}),
and using the tilt formula~\cite{weeks:tilt,Ush,FriPe} we establish
point~(\ref{main:canonical:item}). 
Next, we use Haken's theory of normal surfaces to prove point~(\ref{main:incompressible:item}), 
and we show points~(\ref{main:genus:item}) and~(\ref{main:turaev:viro:item})
by direct arguments.

\paragraph{Triangulations}
Let $N$ be a compact manifold with boundary and let $\calT$ be an ideal
triangulation of $N$. We associate to $\calT$ the graph $\Gamma_{\calT}$ 
whose vertices are the components of $\partial N$ and whose edges 
correspond to the edges of $\cal T$. 
Of course, if $N$ is connected then 
$\Gamma_{\calT}$ is also connected. 

\begin{lemma}\label{graph:lem}
Let $N$ be a connected compact manifold with boundary and let $\calT$ be an ideal
triangulation of $N$. Then $\chi(\Gamma_{\calT})\leqslant 0$. If $\chi(\Gamma_{\calT})=0$
then 

\begin{center}
\begin{tabular}{cc}
$\Gamma_{\calT} = $ & \mettifig{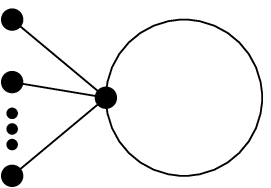, width = 0.8 cm}
\end{tabular}
\end{center}

\noindent and each tetrahedron of $\cal T$
has at least three vertices on the component $C$ of $\partial N$ having
multiple adjacencies in $\Gamma_{\calT}$.
\end{lemma}
\begin{proof}
Each tetrahedron $\tetra$ determines a subgraph $\Gamma_{\tetra}$
of $\Gamma_{\calT}$ whose vertices and edges correspond to the vertices and edges
of $\tetra$, where $\tetra$ is considered as a subset of $N$. 

Now suppose $\chi(\Gamma_{\calT})\geqslant 0$.
Then $\chi(\Gamma_{\tetra})\geqslant 0$
for every $\tetra\in\calT$, and this implies that $\chi(\Gamma_{\tetra})$ is either
{\begin{picture}(21,8)
\put(4.2,4){\circle{8}} \put(8,4){\circle*{2}} \put(8,4){\line(1,0){8}}
\put(16,4){\circle*{2}} \end{picture}}
or
{\begin{picture}(13,8)
\put(4.2,4){\circle{8}} \put(8,4){\circle*{2}}
\end{picture}}.
Therefore each $\tetra\in\calT$
has at least  three vertices on the same
component $C$ of $\partial N$. 
Moreover $\Gamma_\calT$ is the union of the $\Gamma_{\tetra}$'s,
so it is as required and the conclusion follows.
\end{proof}

The following result implies Proposition~\ref{first:prop}. The \emph{incidence number} of an edge
in a triangulation is the number of tetrahedra incident to it (with multiplicity).

\begin{prop}\label{forma:tria:prop}
If $M$ is connected and 
$\partial M = \Sigma_g\sqcup\big(\mathop{\sqcup}\limits_{i=1}^k T_i\big)$ then 
any ideal triangulation $\calT$ of $M$ has at least $g+k$ tetrahedra, and if it has
$g+k$ the following holds:
\begin{itemize}
\item $g\geqslant k$;
\item For any $i=1,\ldots,k$ there are exactly two tetrahedra
of $\calT$ with $3$ vertices on $\Sigma_g$ and one on $T_i$;
the remaining $g-k$ tetrahedra have all $4$ vertices on $\Sigma_g$;
\item
$\calT$ has $k+1$ edges $e_0,\ldots,e_k$ such that $e_0$ 
has both its endpoints on $\Sigma_g$ and incidence number $6g$,
while $e_i$ connects
$\Sigma_g$ to $T_i$ and has incidence number $6$
for any $i=1,\ldots,k$.
\end{itemize}
\end{prop}
\begin{proof}
If $\calT$ has $n$ tetrahedra, an Euler characteristic argument shows that it has
$n-g+1$ edges. Therefore $\chi(\Gamma_{\calT}) = 1+k-(n-g+1) = k+g-n$. 
Lemma~\ref{graph:lem} then implies that $n\geq k+g$, and that if $n=k+g$ 
there exists a component $C$ of $\partial M$ such that 
every tetrahedron has at least $3$ vertices on $C$.
Moreover $\calT$ has $k+1$ edges.
Let $y$ be the number of tetrahedra of $\calT$ having some
(and then one) vertex on $\partial M\setminus C$.

We first claim that $C=\Sigma_g$. Note that $\calT$  
induces on $\partial M\setminus C$ a triangulation
with $k$ vertices and $y$ triangles. If $C\neq\Sigma_g$, 
we would have $2-2g=\chi(\partial M\setminus C)=k-y/2$, whence
$y+4 = 4g + 2k = 2 (g+k) +2g \geqslant y+(g+k)+2g$ and $4\geqslant
3g+k$. Since $g\geqslant 2$, this is a contradiction and our claim is proved.

Having shown that $C=\Sigma_g$, we get $0=\chi(\partial M\setminus C)=k-y/2$, 
so $y=2k$. Therefore the triangulation of $\partial M$ 
induced by $\calT$ has exactly one vertex and two triangles 
on each $T_i$.
So for any $i=1,\ldots,k$ two tetrahedra of $\calT$ 
have one vertex in $T_i$. These $2k$ tetrahedra are distinct, so
$n=g+k\geq 2k$, whence $g\geq k$.
Moreover there is only one edge $e_i$ of $\calT$
incident to $T_i$, both tetrahedra incident to $T_i$
are triply incident to $e_i$, and no other tetrahedron is incident to $e_i$.
So all $e_i$'s have incidence number $6$, and the other edge $e_0$ has 
incidence number $6(g+k)-6k = 6g$,
because a tetrahedron has 6 edges.
\end{proof}

We now turn to the proof of Theorem~\ref{main:teo}.

\paragraph{Geometric tetrahedra}
We prove here Theorem~\ref{main:teo}-(\ref{main:hyperbolic:item}).
In order to construct
a hyperbolic structure on our manifold $M\in\calM_{g,k}$ we realize the tetrahedra of an ideal 
triangulation of $M$ as special geometric blocks in $\matH^3$ and then 
we require that the structures match under the gluings. To describe the blocks
to be used we need some definitions. 

A \emph{partially truncated tetrahedron} is a pair $(\tetra,\calI)$, where 
$\tetra$ is a tetrahedron and $\calI$ is a set of 
vertices of $\tetra$, which are called
\emph{ideal vertices}. In the sequel we will always
refer to $\tetra$ itself as a partially truncated tetrahedron, tacitly implying
that $\calI$ is also fixed. The \emph{topological realization} $\tetra^\ast$ 
of  $\tetra$
is obtained by removing from $\tetra$ the ideal vertices and small open
stars of the non-ideal ones. We call \emph{lateral hexagon} and \emph{truncation 
triangle} the intersection of $\tetra^\ast$ respectively with a face of $\tetra$ 
and with the link in $\tetra$ of a non-ideal vertex. The edges of the 
truncation triangles, which also belong to the lateral hexagons, are called 
\emph{boundary edges}, and the other edges of $\tetra^\ast$ are called
\emph{internal edges}. Note that, if $\tetra$ has ideal vertices, a lateral
hexagon of $\tetra^\ast$ may not quite be a hexagon, because some of its
(closed) boundary edges may be missing. A \emph{geometric realization} of $\tetra$ is an 
embedding of $\tetra^\ast$ in $\matH^3$ such that the truncation triangles
are geodesic triangles, the lateral hexagons are 
geodesic polygons with ideal vertices
corresponding to missing edges, and the truncation 
triangles and lateral hexagons 
lie at right angles to each other.
The classification of the geometric realizations of partially truncated
tetrahedra given in~\cite{FriPe} implies the following facts:
\begin{itemize}
\item
Let $\tetra$ be a partially truncated tetrahedron with one ideal
vertex $v_0$, and take $\alpha\in\matR$ with $0<\alpha<\pi/3$.
Then there exists, up to isometry,
exactly one geometric realization of $\tetra$ 
with dihedral angles $\pi/3$ along the internal edges emanating
from $v_0$, and angle $\alpha$ along the other internal edges; 
this geometric partially truncated tetrahedron will be denoted by 
$\tetra^{\mathrm{id}}_{\alpha}$ (where ``id'' stands for ``ideal'');
\item
Let $\tetra$ be a partially truncated tetrahedron without ideal vertices
and take $\alpha\in\matR$ with $0<\alpha<\pi/3$.
Then there exists, up to isometry, exactly one geometric realization
of $\tetra$ with all the dihedral angles along the internal edges equal
to $\alpha$; this geometric truncated tetrahedron will be denoted 
by $\tetra^{\mathrm{reg}}_{\alpha}$ (where ``reg'' stands for ``regular'').
\end{itemize}

\paragraph{Consistency}
Let $M$ be our manifold in $\calM_{g,k}$,
and let $\calT$ be an ideal triangulation
of $M$ with $g+k$ tetrahedra. We try to give $M$ a hyperbolic structure
with geodesic boundary by realizing the tetrahedra of $\calT$ as 
copies of the geometric polyhedra just described. More precisely,
we know from Proposition~\ref{forma:tria:prop} 
that $\calT$ consists of $2k$ tetrahedra with one vertex on
the boundary tori and $g-k$ tetrahedra 
with all the vertices on $\Sigma_g$.
So we fix $\alpha,\beta\in(0,\pi/3)$, we realize the tetrahedra 
incident to the boundary tori of $M$ as $2k$ copies of
$\tetra^{\mathrm{id}}_{\alpha}$ and the tetrahedra incident only to 
$\Sigma_g$ as $g-k$ copies of $\tetra^{\mathrm{reg}}_{\beta}$.

It was shown in~\cite{FriPe}  that the hyperbolic structure given on
the tetrahedra of $\calT$ extends to the whole of $M$ if and only
if all the matching boundary edges have the same length and the total dihedral
angle around each internal edge is $2\pi$.
Suppose first that $1\leqslant k\leqslant g-1$. In this case
the length condition translates into the equation
$f(\alpha,\beta)=0$, where
$$f(\alpha,\beta)=\frac{\cos^2 \alpha + 1/2}{\sin^2 \alpha}-
\frac{\cos^2 \beta + \cos \beta}{\sin^2 \beta},$$
while, by Proposition~\ref{forma:tria:prop},
the total dihedral angle condition gives the
equation
$$6k\cdot\alpha+6(g-k)\cdot\beta=2\pi.$$
Now let
$\beta(\alpha)=\frac{\pi-3k\cdot\alpha}{3(g-k)}$ be the solution of this equation.
Setting  
$\phi(\alpha)=f(\alpha,\beta(\alpha))$
we easily get that $\lim\limits_{\alpha\to 0} \phi(\alpha)=+\infty$,
$\lim\limits_{\alpha\to \pi/3k} \phi(\alpha)=-\infty$.
Moreover, $\phi$ is strictly monotonic on $(0,\pi/3k)$
so the length and total angle equations have a unique 
solution $(\overline{\alpha}_{g,k},\overline{\beta}_{g,k})$ in $(0,\pi/3)\times (0,\pi/3)$.
This solution determines
a hyperbolic structure with geodesic boundary on $M$.

When $k=0$ or $k=g$ the situation is even simpler, and the scheme just described
easily extends. More precisely, when $k=0$ only compact geometric
polyhedra arise, so the shape of the tetrahedra of $\calT$ is parametrized 
by $\beta$ and the hyperbolicity condition is verified for
$\overline{\beta}_{g,0}=\pi/3g$. 
On the other hand, when $k=g$ every tetrahedron
of $\calT$ has a vertex on a boundary torus, so the geometric realizations
of $\calT$ are parametrized by $\alpha$ and the hyperbolicity condition
gives $\overline{\alpha}_{g,g}=\pi/3g$.

\paragraph{Completeness}
To check completeness of the hyperbolic structure just described
we have to determine the similarity structure it induces on the boundary tori. 
By construction, each torus in $\partial M$ is tiled by two equilateral 
Euclidean triangles. This shows that the structures on the
boundary tori are indeed Euclidean, so the hyperbolic structure constructed
in the previous paragraph is complete, and corresponds by Mostow's
rigidity theorem to the \emph{unique} complete finite-volume hyperbolic structure
with geodesic boundary on the topological manifold $M$ with the boundary tori removed.
The volume of $M$ is $2k\cdot {\rm Vol}(\tetra^{\mathrm{id}}_{\overline{\alpha}_{g,k}})
+ (g-k)\cdot {\rm Vol}(\Delta^{\mathrm{reg}}_{\overline{\beta}_{g,k}})$,
which depends on $g$ and $k$ only. We have eventually proved 
Theorem~\ref{main:teo}-(\ref{main:hyperbolic:item}).

\paragraph{Canonical decomposition}
We now establish Theorem~\ref{main:teo}-(\ref{main:canonical:item}).
Kojima proved in~\cite{Kojima} that a complete finite-volume hyperbolic manifold
$M$ with non-empty geodesic boundary admits a \emph{canonical decomposition} into
partially truncated polyhedra (an obvious generalization of a partially truncated
tetrahedron). This decomposition is obtained
by projecting first to $\matH^3$ and then to $M$ the faces of the convex hull
of a certain family $\calP$ of points in Minkowsky 4-space.
This family $\calP$ splits as $\calP'\sqcup\calP''$, with $\calP'$
consisting of the points on the hyperboloid
$\|x\|^2=+1$ which are dual to the hyperplanes giving $\partial\widetilde{M}$, 
where $\widetilde{M}\subset\matH^3$ is a universal cover of $M$. 
The points in $\calP''$ lie on the light-cone,
and they are the duals of horoballs projecting 
in $M$ to Margulis neighbourhoods of the cusps.
The choice of these Margulis neighbourhoods is somewhat tricky, and carefully
explained in~\cite{FriPe}. It will sufficient for our present purposes to know
that any choice of sufficiently small Margulis neighbourhoods leads to a
set $\calP''$ which works. Note in particular that, when there is more 
than one cusp, the Margulis neighbourhoods need not have the same volume,
as required for instance for the canonical Epstein-Penner decomposition~\cite{EpPe}.
In the sequel we will denote by $\calO$ the union of 
sufficiently small Margulis neighbourhoods of the cusps.

\paragraph{Tilts}
Suppose a geometric triangulation $\calT$ of 
$M$ is given. 
The matter of deciding if $\calT$ is the canonical Kojima decomposition
of $M$ is faced using the \emph{tilt formula}~\cite{weeks:tilt, Ush, FriPe},
that we now briefly describe.

Let $\sigma$ be a $d$-simplex in $\calT$ and $\tilde{\sigma}$ be
a lifting of $\sigma$ to $\widetilde{M}\subset\matH^3$. To each end of $\tilde{\sigma}$ 
there corresponds (depending on the nature of the end) 
one horoball in the lifting of $\calO$ 
or one component of the geodesic boundary of $M$, so $\tilde{\sigma}$
determines $d+1$ points
of $\calP$. Now let two tetrahedra $\tetra_1$ and $\tetra_2$ share a $2$-face
$F$, and let $\widetilde{\tetra}_1,\widetilde{\tetra}_2$ and $\widetilde F$ be 
liftings of $\tetra_1,\tetra_2$ 
and $F$ to $\widetilde{M}\subset\matH^3$ 
such that $\widetilde{\tetra}_1\cap\widetilde{\tetra}_2=\widetilde{F}$.
Let $\overline{F}$ be
the $2$-subspace in Minkowsky 4-space that
contains the three points of $\calP$ determined by $\widetilde{F}$. For $i=1,2$ let 
$\overline{\tetra}^{(F)}_i$ be the half-$3$-subspace bounded by 
$\overline{F}$ and containing the fourth point of $\calP$ determined
by $\widetilde{\tetra}_i$. Then one can show that $\calT$ is canonical if and only if, 
whatever $F,\tetra_1,\tetra_2$, the following holds:
\begin{itemize}
\item The half-$3$-subspaces
$\overline{\tetra}^{(F)}_1$ and $\overline{\tetra}^{(F)}_2$
lie on  distinct $3$-subspaces and their convex hull 
does not contain the origin of Minkowsky $4$-space.
\end{itemize}
The tilt formula computes 
a real number $t(\tetra,F)$ describing the
``slope'' of $\overline{\tetra}^{(F)}$. More precisely,
one can translate the condition just stated 
into the inequality 
$$t(\tetra_1,F)+t(\tetra_2,F)< 0.$$

Coming to the manifolds we are interested in,
let $M\in\calM_{g,k}$, let $\calT$ be a geometric triangulation
of $M$ by $g+k$ partially truncated tetrahedra 
as described above and let $\calO$ be a suitable neighbourhood
of the cusps of $M$. It was shown in~\cite{FriPe} that $\calO$
determines a real number $r_\Delta(v)>0$ for any ideal vertex $v$
of any tetrahedron $\Delta$ in $\calT$.
This number $r_\Delta(v)$ represents the ``height'' 
of the trace in $\Delta$ near $v$ of $\partial\calO$
(except that $r_\Delta(v)\ll 1$ means that $\partial\calO$ is ``very'' high).

Recall now that in $\calT$ all the tetrahedra having vertices in the cusps
have only one such vertex and are isometric to each other.
It follows that we can choose
$\calO$ so that $r_\Delta(v)$ has a certain
value $r$ whenever $v$ is an ideal vertex of some $\Delta$ in $\calT$.
Using the formulae given in~\cite{FriPe} 
we can now easily compute the tilts of the geometric blocks
of $\calT$.
\begin{itemize}
\item
Let $v_0$ be the ideal vertex of $\tetra^{\mathrm{id}}_{\alpha}$ and
let $r=r(v_0)$ be the parameter associated to the intersection of  
$\tetra^{\mathrm{id}}_{\alpha}$ with $\calO$. If $F_0$
is the face of $\tetra^{\mathrm{id}}_{\alpha}$ opposite
to $v_0$ and $F_1$ is any other face of
$\tetra^{\mathrm{id}}_{\alpha}$, then 
\begin{eqnarray*}
t(\tetra^{\mathrm{id}}_{\alpha},F_0)&=&r/(2\cos\alpha)-
\sqrt{4\cos^2\alpha - 1},\\
t(\tetra^{\mathrm{id}}_{\alpha},F_1)&=&-r/2;
\end{eqnarray*}
\item
If $F$ is any face of $\tetra^{\mathrm{reg}}_{\beta}$, then
$$t(\tetra^{\mathrm{reg}}_{\beta},F)=-\sqrt{\frac{(3\cos\beta 
-1)(2\cos\beta - 1)}{\cos \beta + 1}}.$$
\end{itemize} 
If $r$ is small enough, we then get $t(\tetra_1,F_1)+t(\tetra_2,F_2)<0$
for any pair $(F_1,F_2)$ of matching faces of $\calT$, and this suffices
to prove that $\calT$ is the Kojima canonical decomposition of $M$.
Therefore there is only one triangulation of $M$ with $g+k$ tetrahedra. We have
proved Theorem~\ref{main:teo}-(\ref{main:canonical:item}).

\begin{figure}
\begin{center}
\mettifig{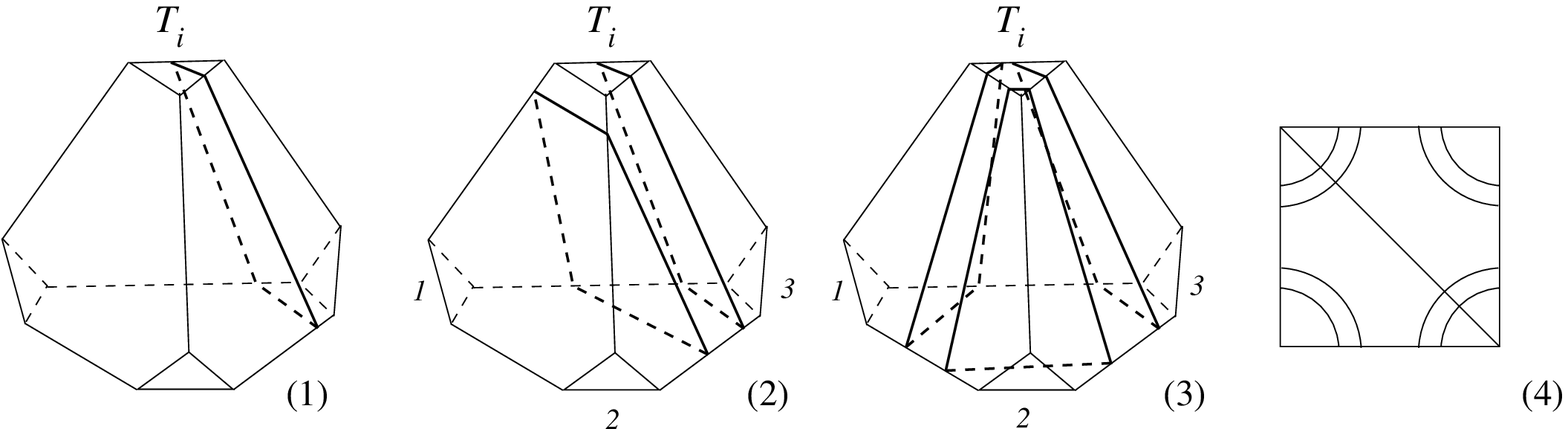, width = 14 cm}
\nota{Types of intersection between a normal surface and a tetrahedron.} 
\label{normal:fig}
\end{center}
\end{figure}

\paragraph{Normal surfaces}
We now prove Theorem~\ref{main:teo}-(\ref{main:incompressible:item}) 
using Haken's theory of normal surfaces~\cite{Mat}. Let $\calT$ be the
ideal triangulation of $M$ with $g+k$ vertices. Suppose $S\subset M$
is a properly embedded incompressible and $\partial$-incompressible 
surface not intersecting $\Sigma_g$. Then $S$ can be isotoped into 
normal position with respect to $\calT$. The surface $S$
intersects each tetrahedron
having $4$ truncation triangles in $\Sigma_g$ in internal triangles and
squares, and it intersects each tetrahedron with $1$ truncation
triangle in some boundary torus $T_i$ in internal triangles,
squares, and squares
having one edge in $T_i$ as in Fig.~\ref{normal:fig}-(1).
We prove that only internal triangles are permitted, which implies that
$S$ is boundary-parallel, hence not essential.

Let us first consider a tetrahedron $\Delta$ with one truncation triangle
in some $T_i$. Suppose first there are $q_1>0$ internal
squares. One type of square intersecting $T_i$ can also be
present, as shown in Fig.~\ref{normal:fig}-(2). Let $q_2$ be the
number of parallel copies of such squares. We enumerate the three
truncation triangles of $\Delta$ not on $T_i$ by 1, 2, and 3 and we denote by $t_j$
the number of triangles in $S\cap\Delta$ that are parallel to the $j$-th
truncation triangle. The three base edges connecting vertices 1,2, and
3 are glued together, therefore we have
$$t_1+t_2 = t_2+q_2+q_1+t_3 = t_1+q_2+q_1+t_3$$
which implies that $t_1=t_2=q_1+q_2+t_3$. The three other edges
are also glued together, whence
$t_1+q_1 = t_3$, a contradiction. This shows that there are
no internal squares. Then $S$ can contain three types of squares intersecting
$T_i$, as shown in Fig.~\ref{normal:fig}-(3), and we denote by $q_1,
q_2,$ and $q_3$ the number of parallel copies of each type. As above,
we have two equations, namely
$$t_1+t_2+q_1+q_2 = t_2+t_3+q_2+q_3 = t_3+t_1+q_3+q_1,$$
$$t_1 = t_2 = t_3,$$
giving $q_1 = q_2 = q_3$. Then $\partial S \cap T_i$
consists of $q_1$ copies of the trivial loop, as in
Fig.~\ref{normal:fig}-(4):
a contradiction, since $S$ is incompressible.

The case where $\Delta$ has all truncation triangles in $\Sigma_g$ is easier:
only triangles and squares are allowed, and all six edges of $\Delta$ are glued
together. Writing equations as above we get that there is no square.

\paragraph{Matveev complexity and Heegaard genus}
Theorem~\ref{main:teo}-(\ref{main:complexity:item}), 
which states that $c(M)=g+k$ for $M\in\calM_{g,k}$, is now an easy
consequence of
Proposition~\ref{first:prop} together with the 
fact~\cite{Mat} that, because of hyperbolicity,
$c(M)$ equals the minimal number of tetrahedra in an ideal
triangulation of $M$.

The genus of $\big(M,\Sigma_g,\mathop{\sqcup}\limits_{i=1}^k T_i\big)$ 
is of course at least
$g$, and it is actually at most $g+1$ because,
if $e$ is the only edge of the minimal triangulation of $M$ having both
ends in $\Sigma_g$, the boundary of a regular neighborhood of 
$\Sigma_g\cup e$ is easily seen to be a Heegaard surface.
The next result, together with 
the fact that $H_{g,k}$ is not hyperbolic,  
shows that the genus is indeed $g+1$.

\begin{lemma}\label{genus:g:lem}
If $M$ is compact with
$\partial M = \Sigma_g\sqcup(\mathop{\sqcup}\limits_{i=1}^k T_i)$ and
$\big(M,\Sigma_g,\mathop{\sqcup}\limits_{i=1}^k T_i\big)$ has genus $g$ then $M=H_{g,k}$.
\end{lemma}

\begin{proof}
$M$ is obtained by attaching $1$-handles to 
$(\mathop{\sqcup}\limits_{i=1}^k T_i)\times[0,1]$ along
$(\mathop{\sqcup}\limits_{i=1}^k T_i)\times\{1\}$ until
a boundary component $\Sigma_g$
is created. Viewing $T_i\times[0,1]$ as the collar of
the boundary of a solid torus, we see that $M$ can also be described
as follows:
\begin{itemize}
\item Attach $1$-handles to a disjoint union of $k$ solid tori until 
a connected manifold with one 
boundary component $\Sigma_g$ is created;
\item Drill the cores of the original $k$ solid tori.
\end{itemize}
At the end of the first step we obviously have $H_g$, so 
we have $H_{g,k}$ at the end of the second step.
\end{proof}

\paragraph{Turaev-Viro invariants}
We conclude this section proving
Theorem~\ref{main:teo}-(\ref{main:turaev:viro:item}). 
As pointed out
in~\cite{Mat:TV}, Turaev-Viro invariants depend only on incidence numbers
between edges and tetrahedra in a triangulation. In our context,
let us consider the minimal triangulation $\calT$ of some $M$ in $\calM_{g,k}$.
If we assign the colour $0$ to the edge having both endpoints in $\Sigma_g$, 
and colours $\{1,\ldots,k\}$ to the other edges,
then the $6$ edges of each of the $g+k$ tetrahedra in $\calT$
are coloured.
It is clear that this set of $g+k$ coloured
tetrahedra is the same for each $M\in\calM_{g,k}$. This implies that
all such $M$'s have the same Turaev-Viro invariants.

\section{Spines and Dehn filling} \label{spines:section}
\begin{figure}
\begin{center}
\mettifig{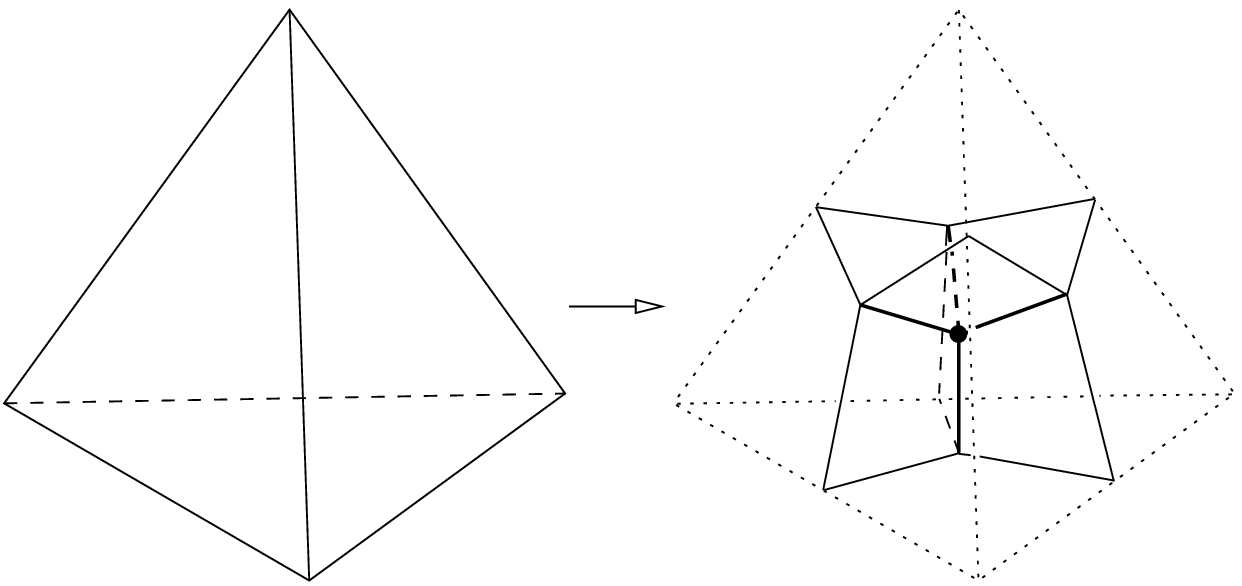, width = 8 cm}
\nota{From an ideal triangulation to a standard spine.} \label{dualspine:fig}
\end{center}
\end{figure}

We prove here Theorem~\ref{main:teo}-(\ref{main:homology:item}),
Theorem~\ref{main:Dehn:teo}, Proposition~\ref{1-bridge:prop}
and Proposition~\ref{tangle:prop}. 
To do this, we switch from the viewpoint of ideal
triangulations to the dual viewpoint of \emph{standard spines},
suggested in Fig.~\ref{dualspine:fig}. Recall that a spine
of a manifold is a subpolyhedron onto which the manifold collapses.
A polyhedron is \emph{standard} if it is locally homeomorphic to that of
Fig.~\ref{dualspine:fig}-right and its natural stratification consists
of $0$-, $1$-, and $2$-cells.
We will be tacitly using in the sequel some of Matveev's 
theory of spines~\cite{Mat}, but we actually will not need to cite 
any precise result: we will try to reconstruct all we need in an elementary
and self-contained way.  

Let us then fix $M\in\calM_{g,k}$ and the spine $P$ dual to the triangulation
of $M$ with $g+k$ tetrahedra. Note that $P$ has 
a cellularization into
vertices, edges, and faces corresponding to tetrahedra, faces, and edges of 
the triangulation. We denote in particular by $S(P)$ the 1-skeleton of $P$ (a 4-valent graph).
By Proposition~\ref{forma:tria:prop} the spine $P$ contains
$k$ (open) hexagonal faces $F_1,\ldots, F_k$ and one big face $G$ with $6g$
vertices (with multiplicity). For $i=1,\ldots,k$ the closure $\overline{F_i}$ of $F_i$ is a torus
which bounds a collar of the $i$-th toric component $T_i$ of $\partial M$,
and the rest of $P$ lies outside this collar.

\paragraph{Homology}
We prove Theorem~\ref{main:teo}-(\ref{main:homology:item}).
The case $k=0$ was dealt with in~\cite{FriMaPe1}, so we suppose $k>0$.
Since $M$ collapses onto $P$,
we have $H_1(M;\matZ)\cong H_1(P;\matZ)$, and
we can use cellular homology to compute $H_1(P;\matZ)$.
Since $\overline {F_i}$ intersects $S(P)$ in a $\theta$-shaped graph,
there is a maximal tree $Y$ in the 4-valent graph $S(P)$ intersecting
each $\overline{F_i}$ in an edge.
Then $S(P)\setminus Y$ consists of $g+k+1$ edges
$e_1,\ldots,e_{g+k+1}$, where $e_{2i-1}$ and $e_{2i}$ 
are contained in $\overline{F_i}$ for $i=1,\ldots,k$, while
$e_{2k+1}, \ldots, e_{g+k+1}$
are contained in $\overline G$.
Choosing an orientation on each $e_j$, $F_i$, and $G$, we get a presentation
for $H_1(P)$ with generators $e_1,\ldots,e_{g+k+1}$ and relators
given by the incidence numbers of $G$ and the $F_i$'s on the $e_j$'s.
Each $F_i$ contributes with
the trivial relator $e_{2i-1}+e_{2i}-e_{2i-1}-e_{2i}$, while $G$ contributes
with a big relator $w$ containing $e_1$ once. Therefore
$H_1(P;\matZ) = \matZ^{g+k+1}/\langle w\rangle\cong\matZ^{g+k}$.

\paragraph{Dehn fillings} 
We prove here Theorem~\ref{main:Dehn:teo}, starting from the case $h=1$
(the general case will easily follow). Let then $\alpha$ be a slope
on the boundary torus $T_1$ corresponding to $\overline {F_1}$.
It is easy to construct a spine $P(\alpha)$ for $M(\alpha)$:
the complement of $P\subset M \subset M(\alpha)$ inside
$M(\alpha)$ consists of the disjoint union of $\Sigma_g\times [0,1)$,
$k-1$ copies of $\Sigma_1\times[0,1)$, and one open solid torus. 
Take a meridinal disc $D$ of this solid torus.
The complement of $P\cup D$ is as above, with an open ball instead of
the open solid torus.
The loop $\partial D\subset \overline{F_1}$ cuts $F_1$ into some open faces of
$P\cup D$ (see two examples in Fig.~\ref{cusp3:fig} below). Each such face
separates $\Sigma_g\times [0,1)$ from the open ball, so, if we remove the face,
we get a spine $P(\alpha)$ of $M(\alpha)$.

The $\theta$-shaped graph $S(P)\cap \overline{F_1}$ contains three
loops, representing three slopes with
pairwise intersection one, having coordinates $0,1,$ and $\infty$
with respect to an appropriate basis of $H_1(T;\matZ)$.
Let us consider the case $\alpha$ is $0,1,$ or $\infty$.
In the construction sketched above of a spine of $M(\alpha)$, 
we can ask $\partial D$ to lie inside the 
$\theta$-shaped graph $S(P)\cap \overline{F_1}$.
The face $F_1$ then survives in $P\cup D$, hence $P(\alpha) = (P\cup D)\setminus
F_1$ is a spine of $M(\alpha)$.
Now $P(\alpha)$ has an induced stratification with 
1-dimensional stratum $S(P)$,
and 2-dimensional stratum consisting of the 
$k-1$ faces $F_2, \ldots, F_k$, the face $G$
and the disc $D$. 
Note now that there is one edge of $P(\alpha)$, which was
previously adjacent twice to $F_1$ and once to $G$, 
which is now only adjacent once to $G$. Therefore $P(\alpha)$ can be
collapsed starting from this edge, and in this collapse the
2-dimensional strata $G$ and $D$ disappear. The resulting polyhedron
is still a spine of $M(\alpha)$, and is made of $k-1$ tori 
$\overline{F_2}, \ldots, \overline{F_k}$ and some 1-dimensional
strata, \emph{i.e.}~a graph connecting these tori. An orientable manifold having
such a spine is necessarily a boundary connected sum of a handlebody and some 
$\Sigma_1\times I$. Since
$\partial M(\alpha)$ consists of one $\Sigma_g$ and $k-1$ tori, we have
$M(\alpha) = H_{g,k-1}$, as required.

We consider now the case $\alpha\in\{-1,1/2,2\}$. 
The slope $\alpha$ is represented by a loop $\partial D$
which intersects transversely the graph $S(P)\cap \overline{F_1}$ in two points,
as in Fig.~\ref{cusp3:fig}-left.
\begin{figure}
\begin{center}
\mettifig{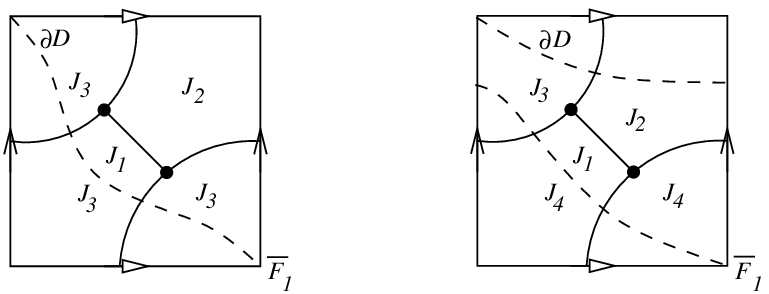}
\nota{A slope in $\{-1,1/2,1\}$ (left) or in
$\{-2,-1/2,1/3,2/3,3/2,3\}$ (right) is represented by a loop $\partial D$
intersecting transversely $S(P)\cap \overline{F_1}$ in two (left) or
three (right) points.} 
\label{cusp3:fig}
\end{center}
\end{figure}
Consider the (open) face
$J_1\subset\overline{F_1}$ shown in Fig.~\ref{cusp3:fig}-left.
The spine $P(\alpha) = (P\cup D)\setminus J_1$ of $M(\alpha)$
has an induced stratification with 
$S(P)\cup\partial D$ as 1-dimensional stratum 
and the faces $F_2,\ldots,F_k$, $G$, $J_2$, $J_3$, and $D$ as 2-dimensional
strata. However this is not the intrinsic stratification of $P(\alpha)$,
because each of the four edges that were adjacent to $J_1$ is now
adjacent to a pair of faces only, so the four edges and the faces incident
to them can be merged into a single 2-dimensional stratum $S$.
The pairs of faces incident to the edges of $J_1$ are
$\{J_3,G\}$, $\{J_2,G\}$, $\{J_3,G\}$, and $\{J_3,D\}$, therefore $S$
is either an annulus or a M\"obius strip.
Let us consider its core $\gamma$. Taking the pre-image of $\gamma$ under
the projection from $M(\alpha)$ to 
$P(\alpha)$ we get an annulus or M\"obius strip $R$ 
properly embedded in $M(\alpha)$ with $\partial R\subset \Sigma_g$
and intersecting $P(\alpha)$ in $\gamma$.
Moreover, cutting $M(\alpha)$ and $P(\alpha)$ along $R$ and $\gamma$ respectively, 
we get a manifold $M'$ and a spine of $M'$
which retracts onto $P(\alpha)\setminus S$.
This polyhedron is easily seen to be connected, so
$M'$ is connected, \emph{i.e.}~$R$ is non-separating.
In addition, $P(\alpha)\setminus S$ consists of 
the tori $\overline{F_2},\ldots, \overline {F_k}$ and some graph connecting them,
which implies as above that $M' = H_{g,k-1}$. 
Consider now the manifold $D(M(\alpha))$ obtained by mirroring
$M(\alpha)$ in $\Sigma_g$, and note that it is hyperbolic if $M(\alpha)$ is.
Now $R$ gives a closed non-separating 
surface $D(R)$ in $D(M(\alpha))$, and $D(R)$ is
homeomorphic either to the torus or to the Klein
bottle. Such a surface cannot exist in a hyperbolic manifold,
so $M(\alpha)$ is not hyperbolic.
 
If $\alpha$ is none of the slopes studied above, then by Thurston's
geometrization theorem either $M(\alpha)$ is hyperbolic or it contains
an essential surface of non-negative Euler characteristic. 
Theorem~\ref{main:teo}-(\ref{main:incompressible:item}), now proved,
implies that $M(\alpha)$ does not contain any closed essential
surface. To conclude we now refer to the bounds on
$\largebound$ stated in Section~\ref{statements:section}.
The fact that the slopes $0,1,\infty$ are of type $D$ and
the bound $\largebound(D,D)\leq 1$ imply that $\alpha$ cannot
be of type $D$. The same fact and the bound
$\largebound(D,A)\leq 2$ imply that $\alpha$ cannot be of type $A$,
and hyperbolicity of $M(\alpha)$ follows. Moreover $M(\alpha)$ 
has genus $g+1$ , because $M$ has genus at most $g+1$,  
but it cannot have genus $g$ by Lemma~\ref{genus:g:lem}.

Finally, suppose $\alpha\in\{-2,-1/2,1/3,2/3,3/2,3\}$.
The slope $\alpha$ is represented by a loop $\partial D$
intersecting transversely the graph $S(P)\cap \overline{F_1}$ in three points,
as in Fig.~\ref{cusp3:fig}-right. As above, we take $P(\alpha)=(P\cup
D)\setminus J_1$. The edges that were adjacent to $J_1$ are now 
adjacent to the four pairs of faces $\{J_3, G\}$, $\{J_2, G\}$, 
$\{J_4,G\},$ and $\{J_4, D\}$. We can therefore as above take a
stratification with a 4-valent graph as 1-stratum and discs
$F_2,\ldots,F_k,D'$ as 2-strata, where $D'$ is the disc obtained by merging
$J_2$, $J_3$, $J_4$, $G$, and the four edges of $J_1$.
Now $P(\alpha)$ is standard, so it can be dualized to an ideal triangulation of
$M(\alpha)$ with $k$ edges and $g+k-1$ tetrahedra. Therefore
$M(\alpha)\in\calM_{g,k-1}$, as required.

The case $h>1$ follows from the case $h=1$, using the fact that $H_{g,k}(\alpha) =
H_{g,k-1}$ for all slopes $\alpha$, and repeating the same argument
above to prove that if
$\alpha_i\not\in\{-1,0,1/2,1,2,\infty\}$ for all $i$ 
then the filled manifold is hyperbolic.
The last assertion of Theorem~\ref{main:Dehn:teo} is a 
direct consequence of 
Theorem~\ref{main:teo}-(\ref{main:incompressible:item}), so the proof of~\ref{main:Dehn:teo}
is now complete.

\begin{rem} \label{flips:rem}
{\em The construction used in the proof of Theorem~\ref{main:Dehn:teo}
to pass from $P$ to $P(\alpha)$ 
can actually be generalized~\cite{MaPe} to any slope $\alpha$. The idea is to
note that $P\setminus F_1$ has a natural $\theta$-shaped ``boundary,''
to take a spine $Q_\alpha$ of the filling solid torus $H_1$ so that $Q_\alpha$
also has a $\theta$-shaped ``boundary'' on $\partial H_1$, 
and the gluing of $\partial H_1$ to $T_1$ (determined by $\alpha$) matches
$\partial Q_\alpha$ to $\partial(P\setminus F_1)$.
If $\alpha\in\{0,1,\infty\}$ the polyhedron $Q_\alpha$ is a meridinal
disc with a longitudinal arc in $\partial H_1$, as in
Fig.~\ref{mobius_triplet:fig}-(1).
If $\alpha\in\{-1,1/2,2\}$ the polyhedron $Q_\alpha$
is the M\"obius triplet shown in Fig.~\ref{mobius_triplet:fig}-(2).
If $\alpha\not\in\{-1,1/2,2\}$, one has to change the $\theta$-shaped
boundary of the M\"obius triplet via some \emph{flips} 
(see Fig.~\ref{mobius_triplet:fig}-(3,4)), each flip
adding a vertex to $Q_\alpha$ as in Fig.~\ref{mobius_triplet:fig}-(5).

\begin{figure}
\begin{center}
\mettifig{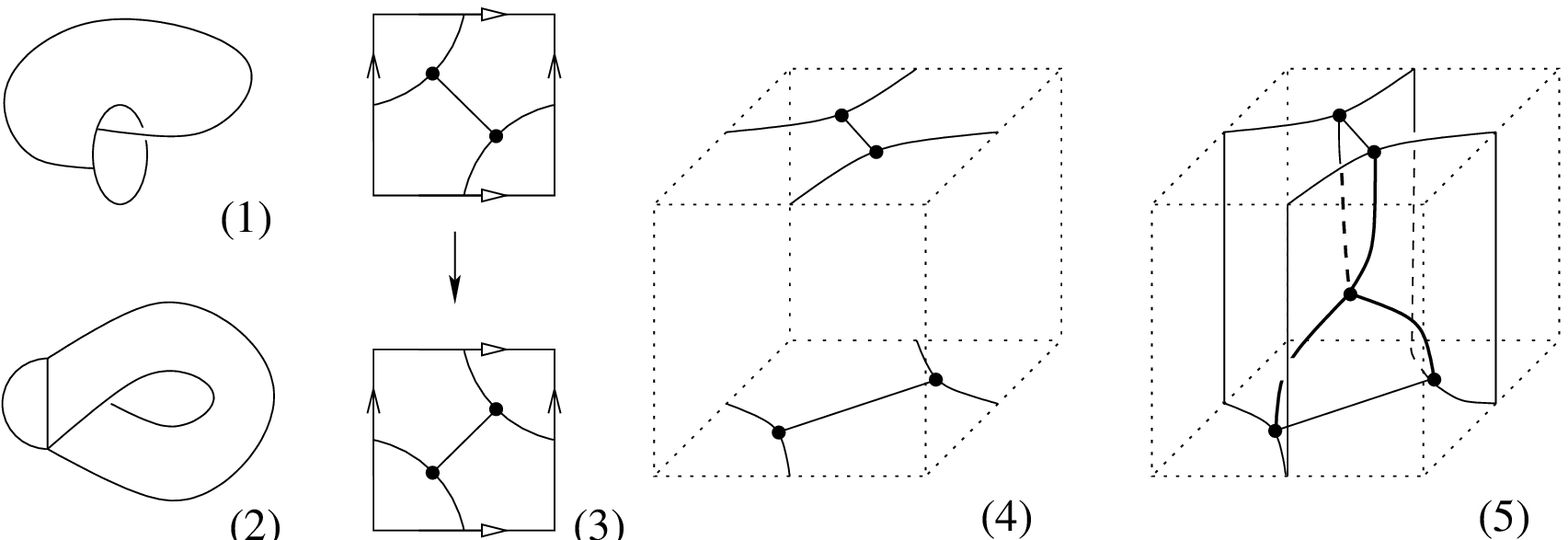, width = 14cm}
\nota{Two spines of the solid torus (1,2) and a flip (3,4) realized by adding a vertex (5).}
\label{mobius_triplet:fig}
\end{center}
\end{figure}

The construction of $P(\alpha)$ as $Q_\alpha\cup(P\setminus F_1)$
is ``efficient,'' in the sense
that if $P$ has a minimal number of vertices, then $P(\alpha)$ very
often does. This construction
is actually dual to adding a layered
solid torus to a triangulation, but it is important to notice that
spines often display greater flexibility than triangulations.
For instance, the construction of $P(\alpha)$
described in the proof when $\alpha\in\{0,1,\infty\}$ has no
analogue for triangulations, is usually efficient, 
and always produces a spine
with strictly fewer vertices than $P$. This is coherent
with the fact that the slopes $\alpha\in\{0,1,\infty\}$ are often exceptional, so they
give a manifold $M(\alpha)$ which is simpler than $M$.
Other natural properties of spines that triangulations do not have are shown in~\cite{Mat}.}
\end{rem}

\paragraph{$1$-bridge knots.}
We now turn to Proposition~\ref{1-bridge:prop}. A knot $K$ in a
manifold $M$ is $1$-bridge
if it can be isotoped to the form $\gamma_0\cup\gamma_1$ where the $\gamma_i$'s
are simple arcs with common ends, $\gamma_0$ lies on
$\partial M$, and $\gamma_1$ is properly embedded and parallel to $\partial M$~\cite{Ga}.

To prove Proposition~\ref{1-bridge:prop}, let $N$ be the exterior
of a knot $K$ contained
in the interior of $H_g$, and assume $N$ is hyperbolic and homeomorphic to some 
$M(\alpha_1,\ldots,\alpha_{k-1})$ for $M\in\calM_{g,k}$. 
Let $P$ be the spine of $M$ dual to
the triangulation with $g+k$ tetrahedra. 
Since $P$ is contained in $M(\alpha_1,\ldots,\alpha_{k-1})$ we can view $P$ as 
as subset of $H_g$. Recall now that $P$ contains
$k$ disjoint hexagonal faces, whose closures are tori,
and one big face $G$. Let $F_1$ be the hexagonal face parallel 
to the only torus in $\partial N$: the torus $\overline{F_1}$ has $K$
on one side and the whole of $P$ on the other side. 
The graph $S(P)\cap\overline{F_1}$ has the 
shape of a $\theta$, so it contains three slopes. Filling along any
of these slopes we get $H_g$, so the bound $\largebound(D,D)\leq 1$ implies,
as in the above proof of Theorem~\ref{main:Dehn:teo}, that the
meridian of $K$ must be one of the slopes contained in $\theta$.

\begin{figure}
\begin{center}
\mettifig{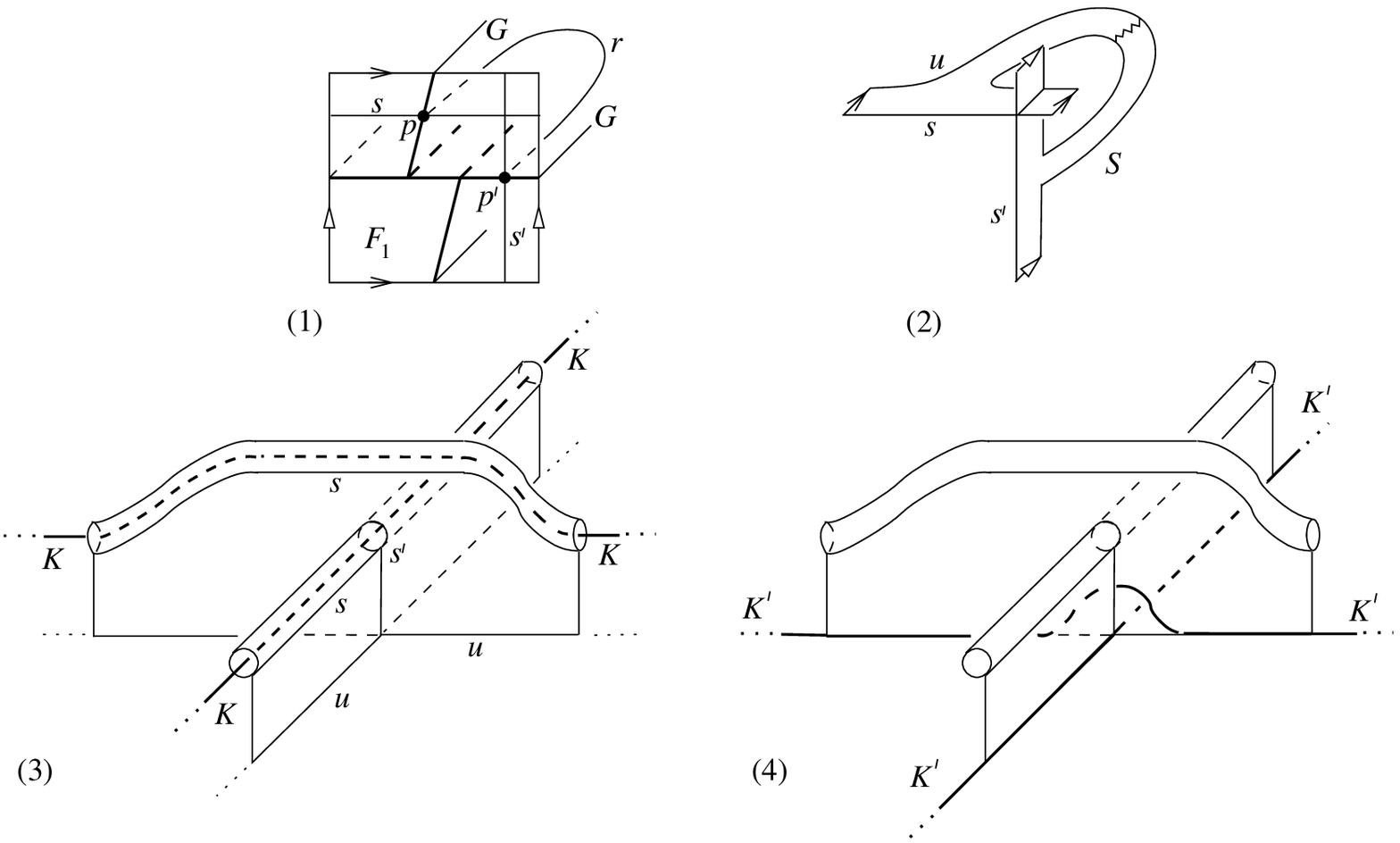, width = 14 cm}
\nota{Isotoping a knot to 1-bridge position.}
\label{1bridge:fig}
\end{center}
\end{figure}

Let us now take in $\overline{F}_1$ two 
loops $s$ and $s'$ so that $s'$ is isotopic to the meridian of $K$ and
$s$ is isotopic to a different slope contained in $\theta$.
We also arrange so that $s$ and $s'$ 
intersect each other and $\theta$ transversely in a single point, 
as in Fig.~\ref{1bridge:fig}-(1).
The points $p=s\cap\theta$ and $p'=s'\cap\theta$ lie in the boundary of the big face $G$,
so there is an arc $r$ properly embedded
in $G$ connecting them, also shown in Fig.~\ref{1bridge:fig}-(1).
The inverse image of the graph $s\cup r\cup s'$ under the retraction of $M$ onto $P$
is a set $S$ as in Fig.~\ref{1bridge:fig}-(2),
where a half-twist of the strip may or not be present along the zig-zagged segment.
In either case, one easily sees that $S$ is a properly 
immersed pair of pants
with $\partial S = s\cup s'\cup u$, where 
$u$ is an immersed loop in $\Sigma_g$ with one self-intersection.
Since $s'$ is a meridian of $K$, the surface $S$ appears in the exterior of $K$ in
$H_g$ as in Fig.~\ref{1bridge:fig}-(3). In part (4) of the same figure we suggest how to isotope $K$
to a knot $K'$ in $1$-bridge position.

\begin{figure}
    \begin{center}
    \mettifig{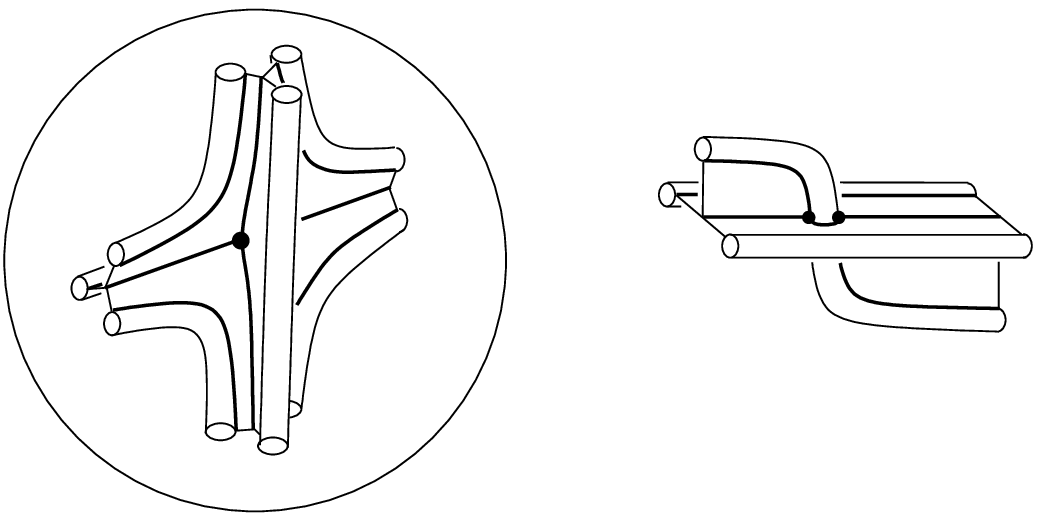, width = 9 cm}
    \nota{Portions of spine of the knot exterior.}
    \label{tangle4:fig}
    \end{center}
    \end{figure}

\paragraph{Knots giving $\calM_{g,1}$}
We prove here Proposition~\ref{tangle:prop}. 
Let $K$ be a knot in $H_g$ constructed with one negative
gluing. We prove that the exterior of $K$ in $H_g$ lies in
$\calM_{g,1}$ by constructing for it a standard spine with $g+1$
vertices. Later we will prove that the unique minimal 
spine of each manifold in
$\calM_{g,1}$ is the result of one such construction, for some $K$.
A spine of the exterior of $K$ is constructed by taking a portion as shown in 
Fig.~\ref{tangle4:fig}-left for each of the $g-1$ tangles. These portions can be
attached (with some torsion) at each positive gluing, while 
the piece shown in Fig.~\ref{tangle4:fig}-right must be inserted 
at the single negative gluing. The resulting polyhedron is
a standard spine of the exterior of $K$ and has $(g-1)+2$ vertices,
as required. Standardness comes from the fact that $K$ is a \emph{knot},
rather than a link, and from the presence of one portion as in 
Fig.~\ref{tangle4:fig}-right.

Now let $M$ be a manifold in $\calM_{g,1}$ and let $P$ be the
spine dual to the triangulation of $M$ with $g+1$ tetrahedra. 
The spine $P$ has one open hexagonal face $F_1$ (whose closure is 
a torus) and one big open face $G$ with $6g$
vertices. The graph $\overline{F_1}\cap S(P)$ has the shape of
a $\theta$, and we denote it by $\theta$. We choose one 
of its edges, say $e$, and distinct points 
$w_1,\ldots,w_{6g-3}$ in the interior of $e$.
We denote by $e_1,\ldots,e_{2g-1}$ the 
edges of $S(P)$ not contained in $\theta$, and choose
an inner point $v_i$ in each $e_i$.
As an abstract face, $G$ is a $6g$-gon with one edge 
$\tilde e$ incident to $e$, two other
edges incident to $\theta$, and $6g-3$ more edges, 
divided into groups of three incident to the same $e_i$.
Let $\tilde w_k$ be the
point of $\tilde e$ incident to $w_k$, and 
$\tilde v_i^{(1)},\tilde v_i^{(2)},\tilde v_i^{(3)}$
be the points on $\partial G$ incident to $v_i$.
Now choose $6g-3$ pairwise disjoint arcs in $G$ 
each having one $\tilde w_k$ and one $\tilde v_i^{(j)}$
as its ends. The image in $P$ of the union
of these arcs is a disjoint union of $2g-1$ graphs,
each having the shape of a $Y$ with all three
endpoints on $e$. Choose now $6g-3$ parallel circles in the
torus $\overline{F_1}$, each intersecting $\theta$ in one of the $w_k$'s.
Attach these circles to the $Y$-shaped graphs, getting $2g-1$ graphs
with shape
\mettifig{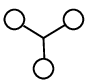, width = .5 cm}. 
It is now not difficult to see that, cutting along these graphs, one gets
$g-1$ polyhedra as in Fig.~\ref{tangle4:fig}-left
and one polyhedron as in Fig.~\ref{tangle4:fig}-right.
(The special portion is the one which contains the two edges of $\theta$ 
other than $e$). From this decomposition one 
readily sees that $P$ arises as explained above for
some knot in $H_g$ constructed with one negative gluing, so $M$
is the exterior of this knot.

\section{Growth estimates}\label{estimates:section}

This section is devoted to the proof of Theorem~\ref{growth:teo}.
Estimates of the form $\#\calM_{g,0}\geq a\cdot b^g$
were already obtained in~\cite{FriMaPe1}. We first improve this
result, and then we extend it to $\calM_{g,k}$ for any fixed $k$.

\begin{teo}\label{factorial:teo}
The sequence $\big(\#\calM_{g,0}\big)_{g=2}^\infty$ has growth type $g^g$.
\end{teo}

Recall that $M\in\calM_{g,0}$ if and only if it is orientable and 
admits a 1-edged
ideal triangulation with $g$ tetrahedra, and that this triangulation is unique.
To prove Theorem~\ref{factorial:teo} we introduce the set $\calG_n$ of homeomorphism classes
of connected $4$-valent graphs with $n$ vertices, and we 
denote by $\calG'_n$ the graphs in $\calG_n$ arising as dual skeleta of
one-edged triangulations of orientable manifolds, 
so that $\#\calM_n\geq\#\calG'_n$. And we prove the following:

\begin{prop}\label{graph:growth:prop}
The sequence $\big(\#\calG_n\big)_{n=1}^\infty$ has growth type $n^n$.
\end{prop}

\begin{prop}\label{phi:map:prop}
For all $n$ there exists a map
$\phi_n:\calG_{n-1}\to\calG'_{n}$
such that $\phi_n(G)$ is obtained by adding a curl at some internal
point of an edge of $G$.
\end{prop}

Assuming these results the proof of Theorem~\ref{factorial:teo} is now easy.
Note that if $G\in\calG_n$ then $G$ contains at most $n$ curls.
So $\#(\phi_n^{-1}(G))\leq n$ for $G\in\calG'_n$. In particular
$\#\calG'_n\geq\frac 1n\cdot\#\calG_{n-1}$, and the conclusion readily
follows from Proposition~\ref{graph:growth:prop} and the next easy:

\begin{rem}\label{max:spines:rem}
\emph{There are at most $18^n$ distinct orientable triangulations with a given
dual $1$-skeleton}.
\end{rem}

\vspace{5pt}\noindent\emph{Proof of~\ref{graph:growth:prop}.}
This result is purely graph-theoretical, and its proof is not hard.
Let us imagine a 4-valent graph with $n$ vertices as being constructed from
the disjoint union of $n$ crosses \ \ 
\begin{picture}(8,8)
\put(2,-4){\line(0,1){12}}
\put(-4,2){\line(1,0){12}}
\put(2,2){\circle*{2}}
\end{picture} \
by joining together in pairs the $4n$ free germs of edges. If we
fix an ordering on these $4n$ germs, there are $4n-1$ choices for the 
germ to be joined to the first germ, then $4n-3$ for the next free germ, and
so on, whence $(4n-1)!!$ in all. This is however too rough a counting, because
disconnected graphs may arise. But the final graph is disconnected if and only
if at some point along the construction process a subgraph without free germs of
edges is created. Assume this happens at time $i$.
For the final edge of the saturated subgraph there is only one choice, so
all other $4n-2i$ choices do not create saturated subgraphs. This easily implies
that at least $(4n-2)!!$ different construction patterns lead to
connected graphs.

We must now consider that different construction patterns can lead to
homeomorphic graphs. Since there are $n$ vertices of valence $4$, one
readily sees that at most $(4!)^n\cdot n!$ different patterns can lead
to the same graph. This implies that
$$(4n-1)!! \geq \#\calG_n\geq \frac{(4n-2)!!}{(4!)^n\cdot n!}.$$
Now the easy inequalities 
$\sqrt{(k+1)!}\geq k!!\geq \sqrt{k!}$ and Stirling's formula imply that
$$\left(\sqrt{8\pi n}\left(\frac{4n}{{\rm e}}\right)^{4n}{\rm e}^{1/48n}\right)^{1/2}
\geq \#\calG_n \geq \frac {\Big(\sqrt{2\pi(4n-2)}\big((4n-2)/{\rm e}\big)^{4n-2}\Big)^{1/2}}
{(4!)^n\sqrt{2\pi n}(n/{\rm e})^n{\rm e}^{1/12n}}$$
and the conclusion readily follows.
\finedimo

To establish Proposition~\ref{phi:map:prop} we begin with the following:

\begin{lemma}\label{1:or:2:lem}
For all $G\in\calG_n$ there exists an orientable ideal triangulation with 
dual graph $G$ and at most two edges.
\end{lemma}

\begin{proof} We adopt the dual viewpoint of orientable standard spines, so we assume
$P$ is a spine such that $S(P)=G$ and the number $m$ of faces of $P$ is the minimal
possible one, and we show that $m\leq 2$.

Recall now that $P$ is determined by a regular neighbourhood $U(S(P))$ of the singular
set $S(P)$, because $\partial U(S(P))$ consists of the attaching circles of the
faces of $P$. Moreover, as in~\cite{BePe}, one can represent $P$ in the plane 
by drawing $\partial U(S(P))$ only.
To establish the conclusion we first prove three claims, each showing
that if $P$ has ``too many'' faces then a new $P$ with fewer faces and the
same $S(P)$ can be constructed, which contradicts minimality.

\emph{Claim 1}: if $e$ is an edge of $S(P)$, the three faces of $P$
running along $e$ cannot be distinct from each other.  This is shown in 
Fig.~\ref{XYZclaim:fig}.
\begin{figure}
\begin{center}
\mettifig{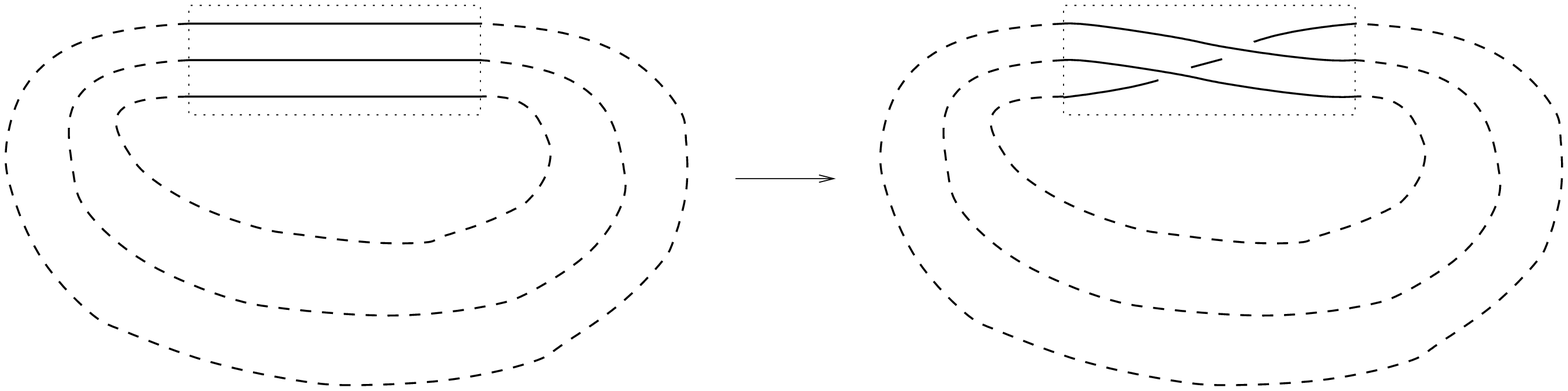,width=12cm}
\nota{Left: if three different faces of $P$ run along the edge $e$
shown in the dotted box, their attaching circles in the 
rest of $\partial U(S(P))$
are connected as shown. Right: this modification of $U(P)$ near $e$
reduces the number of faces.} \label{XYZclaim:fig}
\end{center}
\end{figure}
\emph{Claim 2:} if two faces run along an edge, the face running
twice cannot run in opposite directions, as proved in Fig.~\ref{XYclaim:fig}.
\begin{figure}
\begin{center}
\mettifig{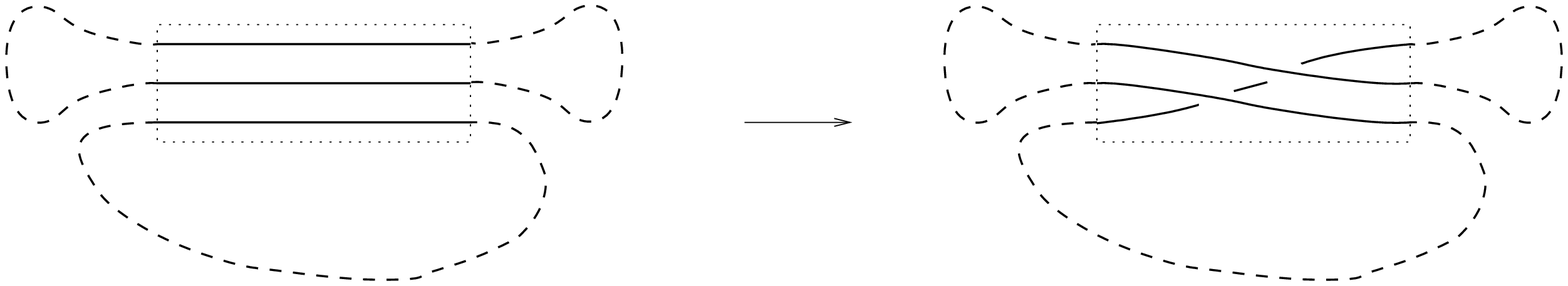,width=12cm}
\nota{If a face passes twice with a counterpass, 
the number of faces can be reduced.}\label{XYclaim:fig}
\end{center}
\end{figure}
\emph{Claim 3:} no face
can run twice along two edges with different companions. This is proved
in Fig.~\ref{XXYZclaim:fig} (note that Claim 2 is used to draw the picture).
\begin{figure}
\begin{center}
\mettifig{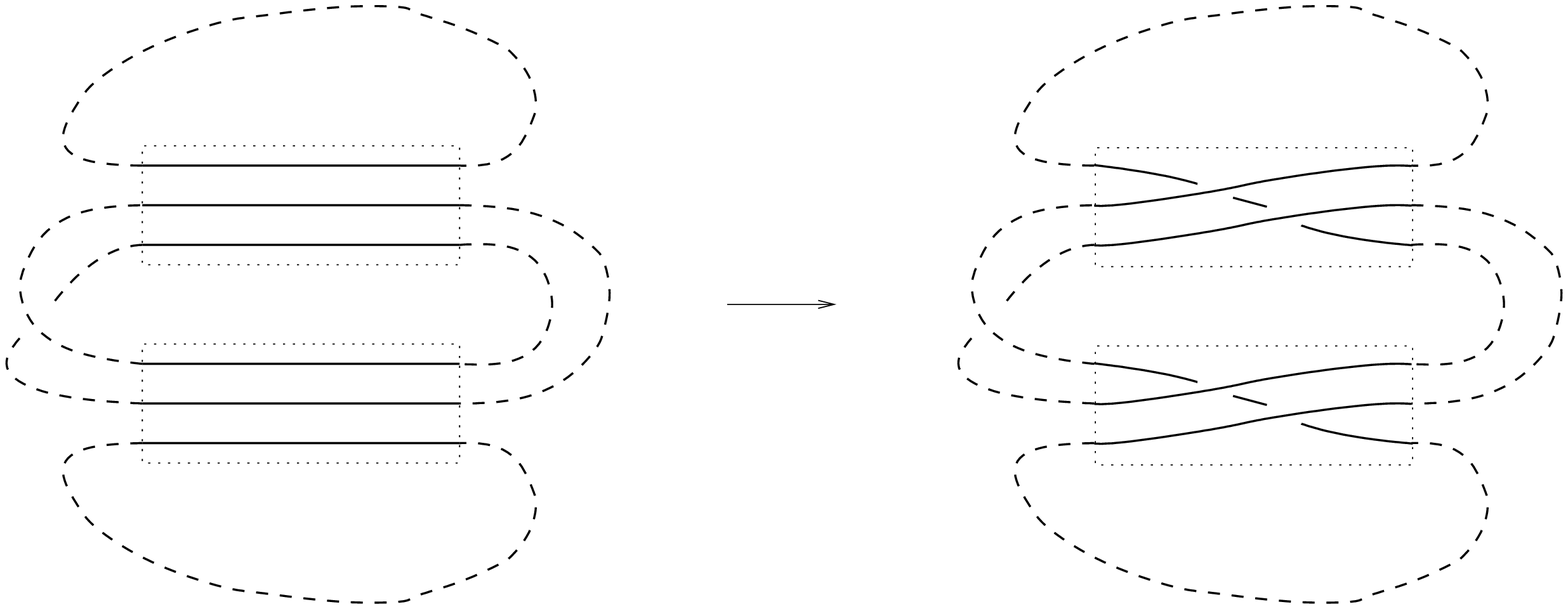,width=12cm}
\nota{If a face passes twice along two edges with different companions,
the number of faces can be reduced.}\label{XXYZclaim:fig}
\end{center}
\end{figure}

We can now conclude. Supposing $P$ has at least three faces, it is easy to see
that there is a simple path $e_1\cdots e_k$ of edges of $S(P)$ such that the total
number of faces touching $e_1\cup e_k$ is at least three. 
Now suppose $k$ is minimal,
let $e_1$ be touched by faces $XXY$ (claim 1 is used here)
and $e_k$ be touched by some $Z\neq X,Y$.
We suppose  $k\geq 3$, leaving the easier case $k=2$ as an exercise
to the reader.
Minimality of $k$ and Claim 3 imply that $e_2$ is touched either by
$XXX$ or by $YYY$, and
actually that the same face $X$ or $Y$ must be triply incident to 
all $e_2,\ldots,e_{k-1}$. Claim 3 now implies that one of the following
must happen:
\begin{enumerate}
\item $e_1=XXY$, $e_2=\ldots=e_{k-1}=XXX$, $e_k=XZZ$;
\item $e_1=XXY$, $e_2=\ldots=e_{k-1}=YYY$, $e_k=YYZ$;
\item $e_1=XXY$, $e_2=\ldots=e_{k-1}=YYY$, $e_k=YZZ$.
\end{enumerate}
Cases (1) and (2) are symmetric, so we treat (2). 
By Claim 2, at the end of $e_1$, the situation
is as in Fig.~\ref{conclusion:fig}-left. 
Since $e_k$ is touched by $YYZ$, Claim 3 shows that the missing label must be
$X$, but then Claim 2 would be violated: a contradiction.
Case (3) is similar. At the end of $e_1$ we still have
\begin{figure}
\begin{center}
\mettifig{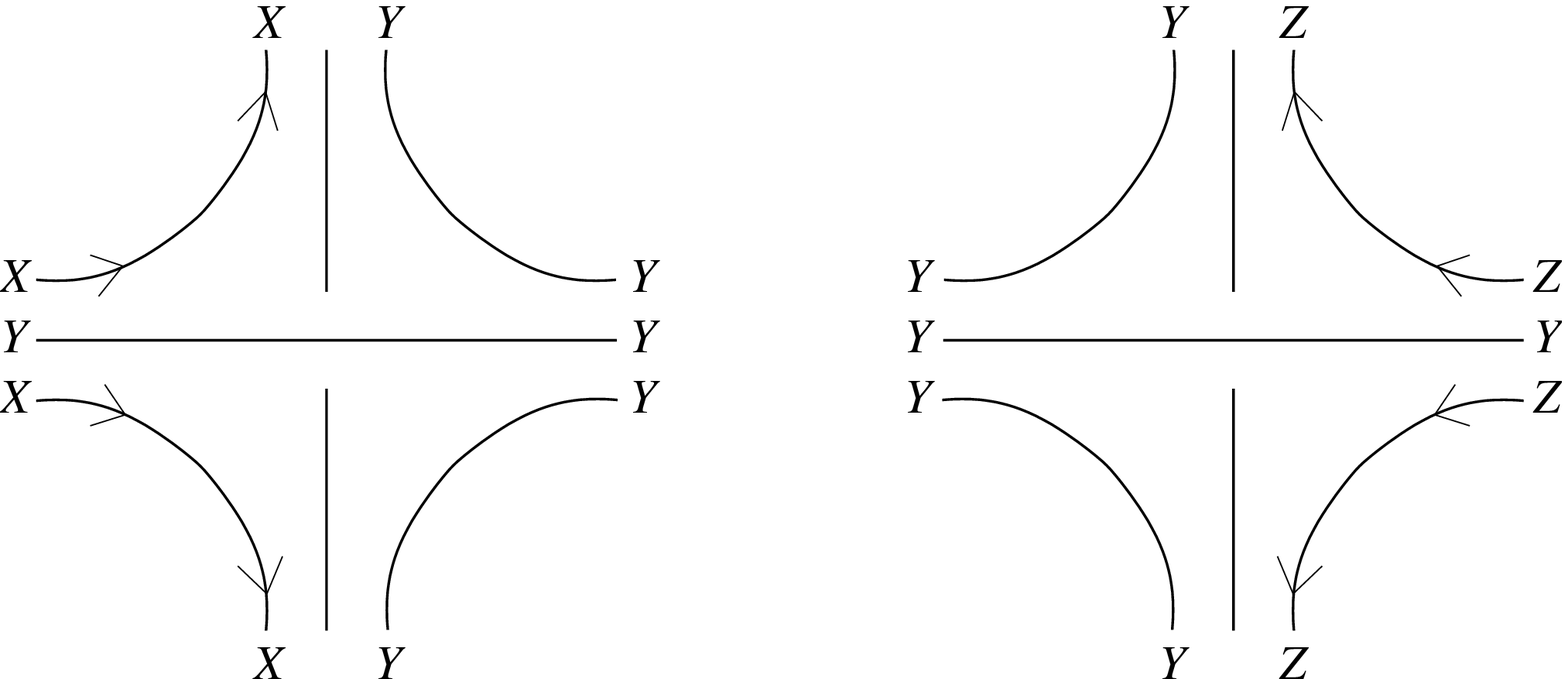,width=8cm}
\nota{Conclusion of the proof.}\label{conclusion:fig}
\end{center}
\end{figure}
the pattern of Fig.~\ref{conclusion:fig}-left, and Claim 2
shows that the missing label must be $Y$. Now the
beginning of $e_k$ is as in Fig.~\ref{conclusion:fig}-right, so
the missing label again must be $Y$. Then we would have edges $XYY$ and $ZYY$, violating
Claim 3.\end{proof}

\begin{figure}
\begin{center}
\mettifig{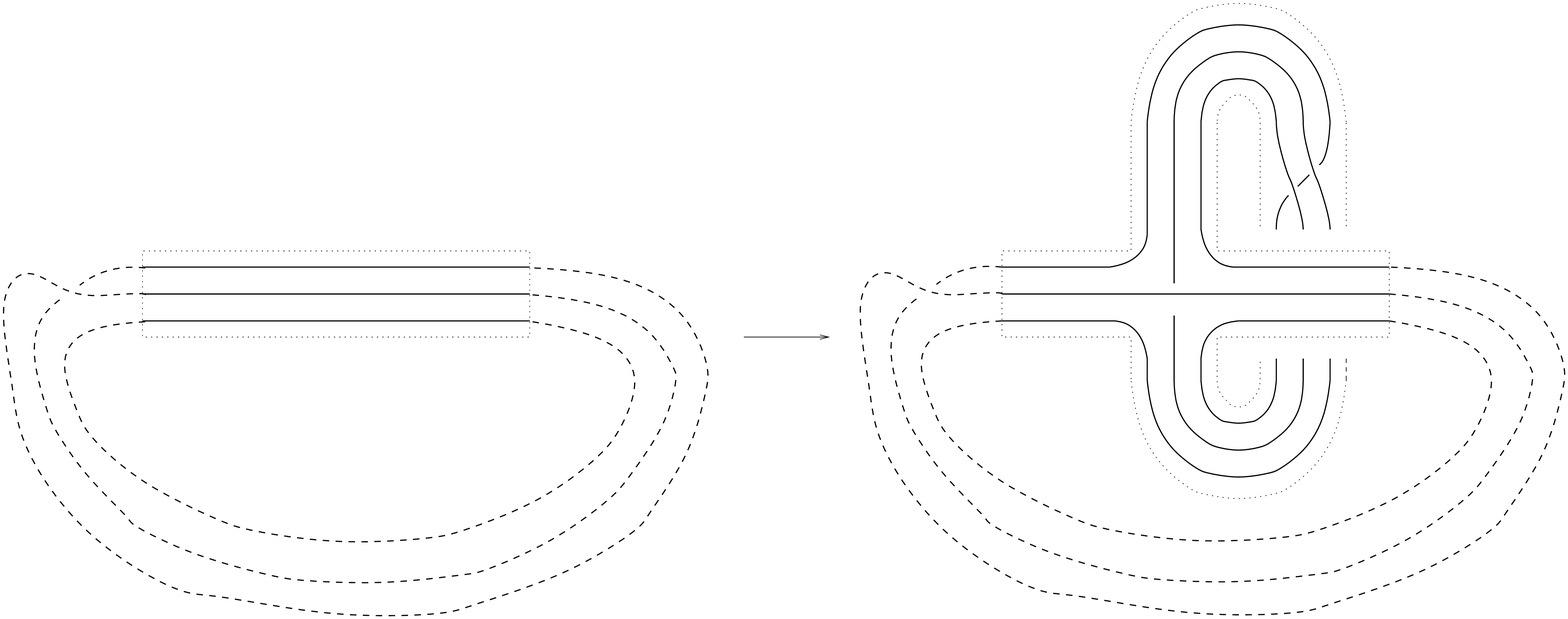,width=14cm}
\nota{Left: attaching circles of faces outside an edge where a face passes twice.
Right: a spine with one face and one more vertex.} \label{curl1:fig}
\end{center}
\end{figure}

\vspace{5pt}\noindent\emph{Proof of~\ref{phi:map:prop}.}
The definition of $\phi(G)$ is different depending on whether
$G$ belongs to $\calG'_{n-1}$ or not. If $G\not\in \calG'_{n-1}$
we consider an orientable standard spine $P$ based on $G$ 
and having two faces. The previous proof implies that there is an edge
of $G$ along which a face of $P$ passes twice in the same direction.
Then we build a new $P$ as in Fig.~\ref{curl1:fig}, and call $\phi_n(G)$ its singular set.
If $G\in\calG'_{n-1}$ we choose $P$ with one face. Looking at any vertex
we see
that there must be an edge along which the face runs twice in one
direction and once in the opposite direction. Then we define $\phi_n(G)$ as in
\begin{figure}
\begin{center}
\mettifig{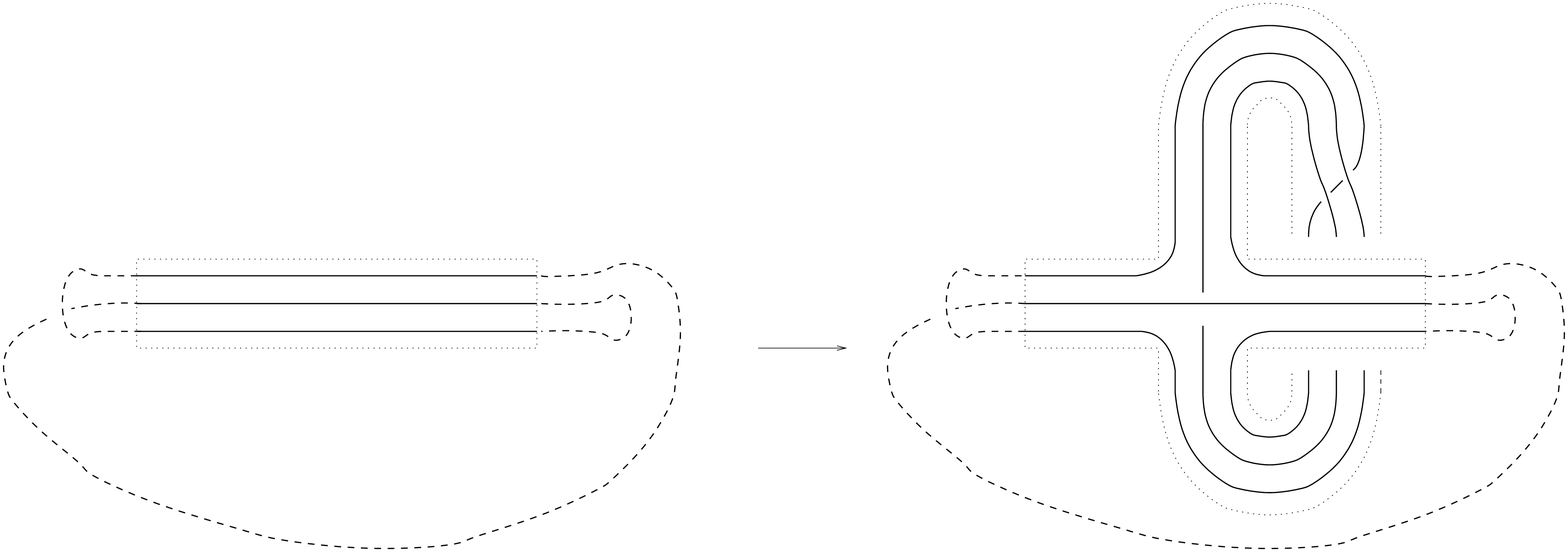,width=14cm}
\nota{Left: attaching circle of the face outside an edge where it passes three
times but not in the same direction.
Right: a spine with one face and one more vertex.} \label{curl2:fig}
\end{center}
\end{figure}
Fig.~\ref{curl2:fig}.\finedimo

\begin{figure}
\begin{center}
\mettifig{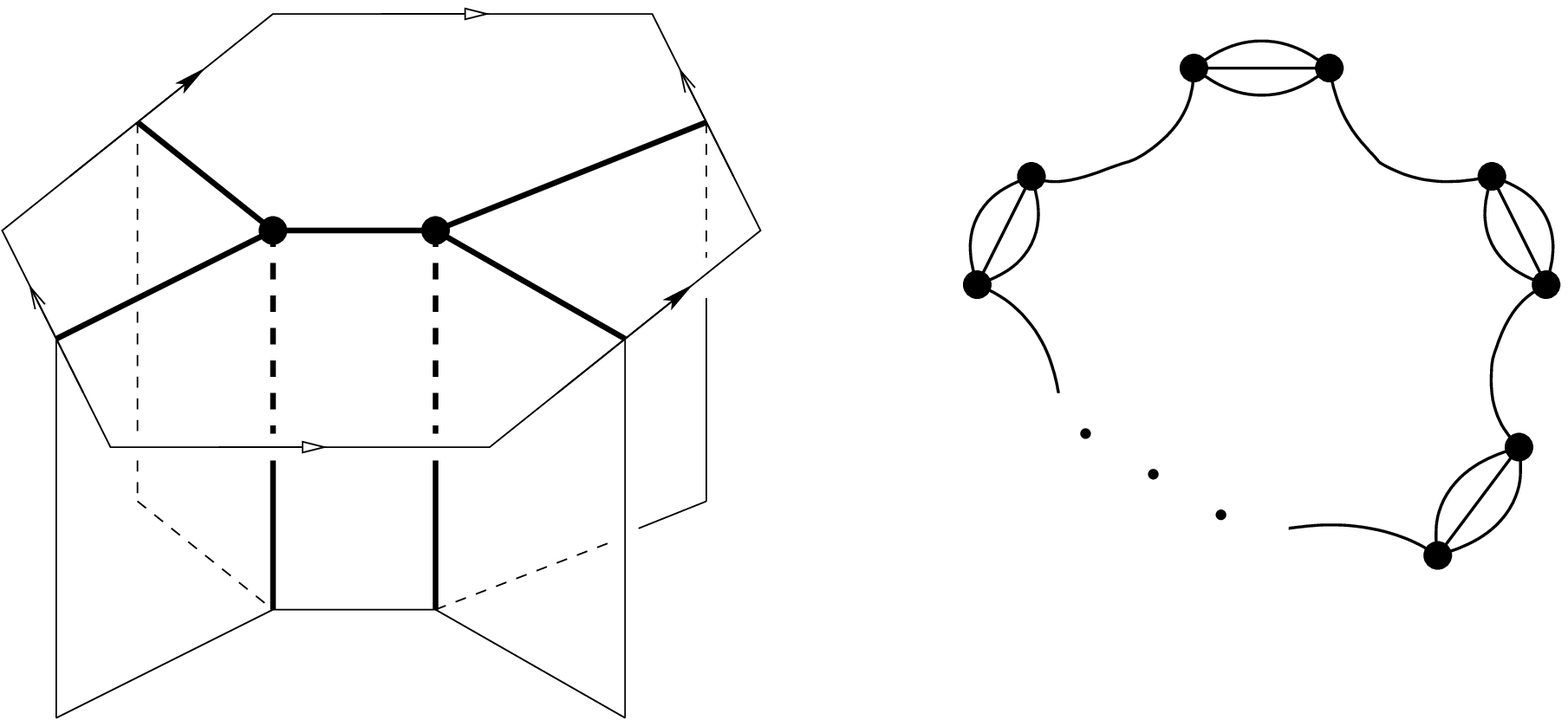,width=8 cm}
\nota{Spines of elements of $\calM_{g,g}$.} \label{new_gequalk:fig}
\end{center}
\end{figure}

\vspace{5pt}\noindent\emph{Proof of~\ref{non-empty:prop} and~\ref{growth:teo}.}
We have already proved in Proposition~\ref{forma:tria:prop} 
that $\calM_{g,k}$ is empty whenever $g<k$.
We also know that if $P$ is a spine with $g+k$ vertices 
of some $M\in\calM_{g,k}$ then up to isotopy $P$ contains the boundary 
tori $T_1,\ldots,T_k$, and a neighbourhood in $P$ of $T_i$ is as in 
Fig.~\ref{new_gequalk:fig}-left.  
Therefore if $g=k$ then $S(P)$ must be as in 
Fig.~\ref{new_gequalk:fig}-right.  This implies that $\calM_{g,g}$ is
non-empty if and only if there exists a spine $P$ of an orientable
manifold having $S(P)$ as in Fig.~\ref{new_gequalk:fig}-right and
a total of $g+1$ faces with $g$ of them as in Fig.~\ref{new_gequalk:fig}-left.  
Using the techniques of~\cite{BePe} it is now easy to see that such a $P$ exists
if and only if $g$ is even.

We now turn to the case $g>k$.
It was proved in~\cite{FriMaPe1} that $\calM_{g,0}$ is non-empty for
all $g\geq 2$, and in~\cite{FriMaPe2}
that $\calM_{2,1}$ is also non-empty.
We will now construct for $g>k$ a function $\psi:\calM_{g,k}\to\calM_{g+1,k+1}$,
whose existence then implies that $\calM_{g,k}$ is non-empty for all $g>k$.
We will also prove that $\#\psi^{-1}(M)\leq 3g$ for all $M$, which,
using Theorem~\ref{factorial:teo}, proves that for
fixed $k$ the growth type
of $\big(\#\calM_{g,k}\big)_{g=k}^\infty$ is $g^g$.

\begin{figure}
\begin{center}
\mettifig{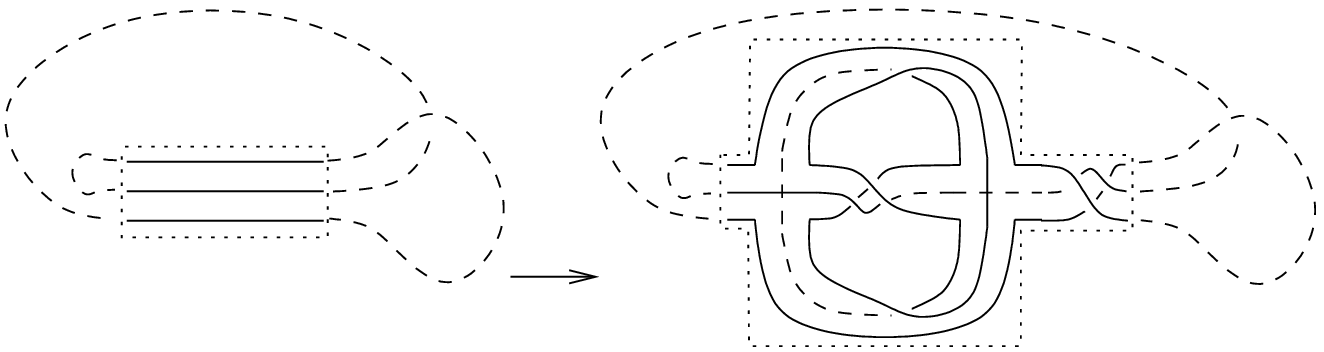, width=12 cm}
\nota{How to transform a spine of $M\in\calM_{g,k}$ to a spine of
$M'\in\calM_{g+1,k+1}$.} \label{add_cusp:fig}
\end{center}
\end{figure}

Let us then construct $\psi$. Let $M\in\calM_{g,k}$ and let 
$P$ be its spine with $g+k$ vertices. Since
$g>k$, there is one vertex $v$ of $P$ which does not belong 
to the closure of any
hexagonal face, so it is adjacent to the big face
$G$ only. Among the 4 edges incident to $v$ there is certainly one edge along
which $G$ runs twice in one direction and once in the opposite one, 
as in Fig.~\ref{add_cusp:fig}-left. With the move shown in
Fig.~\ref{add_cusp:fig} we get a polyhedron $P'$ with one big face,
$g+1$ hexagons and $g+k+2$ vertices. Such a polyhedron is then a spine
of a manifold $M'=\psi(M)\in\calM_{g+1,k+1}$. 
To show that $\#\psi^{-1}(M')\leq 3g$ we recall that
$M'$ has a unique minimal spine $P'$ and note that there are $k+1\leq g$
$\theta$-shaped portions of $P'$ that we can delete and, after deletion,
we have at most $3$ ways to connect what is left.
\finedimo


\begin{thebibliography}{99}


\bibitem {BePe} \textsc{R.~Benedetti -- C.~Petronio}, \emph{A finite
graphic calculus for $3$-manifolds}, Manuscripta Math.~\textbf{88} (1995),
291-310.

\bibitem {Be} \textsc{J.~Berge},
\emph{The knots in $D^2\times S^1$ which have nontrivial Dehn surgeries
that yield $D^2\times S^1$},
Topology Appl.~\textbf{39} (1991), 1--19.

\bibitem {FriMaPe1} \textsc{R.~Frigerio -- B.~Martelli -- C.~Petronio},
\emph{Complexity and Heegaard genus of an infinite class of hyperbolic $3$-manifolds},
{\tt math.GT/0206156}, to appear in Pacific J. Math.

\bibitem {FriMaPe2} \textsc{R.~Frigerio -- B.~Martelli -- C.~Petronio},
\emph{Small hyperbolic $3$-manifolds with geodesic boundary},
{\tt math.GT/0211425}.

\bibitem {FriPe} \textsc{R.~Frigerio -- C.~Petronio},
\emph{Construction and recognition of hyperbolic $3$-manifolds
with geodesic boundary},  {\tt Math.GT/0109012}.

\bibitem{EpPe}\textsc{D.~B.~A.~Epstein -- R.~C.~Penner},
\textit{Euclidean decompositions of noncompact hyperbolic manifolds},
J. Differential Geom. \textbf{27} (1988), 67--80.
                                                    
\bibitem {Ga2} \textsc{D.~Gabai},
\emph{Surgery on knots in solid tori},
Topology \textbf{28} (1989), 1--6.

\bibitem {Ga} \textsc{D.~Gabai},
\emph{$1$-bridge braids in solid tori},
Topology Appl.~\textbf{37} (1990), 221--235.

\bibitem {Go} \textsc{C.~McA.~Gordon}, 
\emph{Small surfaces and Dehn filling}, in
``Proceedings of the Kirbyfest'' (Berkeley, CA, 1998), pp. 177-199 (electronic), 
Geometry and Topology Monographs, Vol. 2, Coventry, 1999. 

\bibitem {GoWu-AA} \textsc{C.~McA.~Gordon -- Y.~Q.~Wu},
\emph{Annular Dehn fillings}, Comment.~Math.~Helv.~\textbf{75} (2000), no. 3, 430--456.

\bibitem {GoWu-DA} \textsc{C.~McA.~Gordon -- Y.~Q.~Wu},
\emph{Annular and boundary reducing Dehn fillings},
Topology \textbf{39} (2000), 531--548.

\bibitem{Kojima}\textsc{S.~Kojima},
\textit{Polyhedral decomposition of hyperbolic manifolds with boundary},
Proc. Work. Pure Math. \textbf{10} (1990), 37-57.

\bibitem {MaPe} \textsc{B.~Martelli -- C.~Petronio}, \emph{$3$-manifolds
up to complexity $9$}, Experiment.~Math. \textbf{10} (2001), 207-236.                             

\bibitem {Mat} \textsc{S.~V.~Matveev}, \emph{Complexity theory of
three-dimensional manifolds}, Acta Appl.~Math.~\textbf{19} (1990),
101-130.

\bibitem {Mat:TV} \textsc{S.~V.~Matveev -- T.~Nowik},
\emph{On $3$-manifolds having the same Turaev-Viro invariants},
Russian J. Math. Phys. \textbf{2} (1994), 317-324.

\bibitem {MiMo} \textsc{K.~Miyazaki -- K.~Motegi},
\emph{Toroidal and annular Dehn surgeries of solid tori},
Topology Appl.~\textbf{93} (1999), 173--178.

\bibitem {TuVi} \textsc{V.~G.~Turaev -- O.~Viro},
\emph{State sum invariants of $3$-manifolds and quantum $6j$-symbols},
Topology \textbf{31} (1992), 865-902.   

\bibitem{Ush}\textsc{A.~Ushijima},
\textit{The tilt formula for generalized simplices in hyperbolic space},
Discrete Comput.~Geom. \textbf{28} (2002), 19-27.

\bibitem{weeks:tilt}\textsc{J.~R.~Weeks},
\textit{Convex hulls and isometries of cusped hyperbolic $3$-manifolds},
Topology Appl. \textbf{52} (1993), 127-149.

\bibitem {Wu2} \textsc{Y.~Q.~Wu},
\emph{Incompressibility of surfaces in surgered $3$-manifolds},
Topology \textbf{31} (1992), 271-279.

\bibitem {Wu} \textsc{Y.~Q.~Wu}, 
\emph{Sutured manifold hierarchies, essential laminations, and Dehn surgery},
J.~Differential Geom.~\textbf{48} (1998), 407--437.

\end{thebibliography}
\end{document}